\documentclass[11pt]{amsart}
\usepackage[english]{babel}
\usepackage{egothic}
\usepackage[T1]{fontenc}

\usepackage{amsmath}
\usepackage{amsthm}
\usepackage{amssymb}
\usepackage[all]{xy}
\usepackage{mathrsfs}
\usepackage{enumerate}
\usepackage{physics}
\usepackage[notcite, final, notref]{showkeys}
\usepackage{dsfont}

\usepackage[dvips]{color}
\newtheorem{theorem}{Theorem}[section]
\newtheorem{problem}[theorem]{Problem}
\newtheorem{lemma}[theorem]{Lemma}
\newtheorem{proposition}[theorem]{Proposition}
\newtheorem{corollary}[theorem]{Corollary}
\newtheorem{definition}[theorem]{Definition}
\usepackage{fdsymbol}
\theoremstyle{definition}
\newtheorem{remark}[theorem]{Remark}

\newtheorem{example}[theorem]{Example}

\newcommand{\ow}{\overline{w}}
\usepackage{xcolor}

\title[Discrete analytic functions and structured matrices]{Discrete analytic functions, structured matrices and  a new family of moment problems}

\author[D. Alpay]{Daniel Alpay}
\address{(DA) Schmid College of Science and Technology \\
Chapman University\\
One University Drive
Orange, California 92866\\
USA}
\email{alpay@chapman.edu}

\author[F. Colombo]{Fabrizio Colombo}
\address{(FC) Politecnico di
Milano\\Dipartimento di Matematica\\Via E. Bonardi, 9\\20133 Milano\\Italy}
\email{fabrizio.colombo@polimi.it}

\author[K. Diki]{Kamal Diki}
\address{(KD) Schmid College of Science and Technology \\
Chapman University\\
One University Drive
Orange, California 92866\\
USA}
\email{diki@chapman.edu}

\author[I. Sabadini]{Irene Sabadini}
\address{(IS) Politecnico di
Milano\\Dipartimento di Matematica\\Via E. Bonardi, 9\\20133 Milano\\Italy}
\email{irene.sabadini@polimi.it}

\author[D.  Volok]{Dan Volok}
\address{(DV) Mathematics Department\\ Kansas State University\\ 138 Cardwell Hall\\1228 N.  17th Street\\
  Manhattan, KS 66506, USA}
\email{danvolok@math.ksu.edu}

\begin{document}
\maketitle

\begin{abstract}
  Using Zeilberger generating function formula for the values of a discrete analytic function in a quadrant we make connections with the theory of
  structured reproducing kernel spaces, structured matrices and a generalized moment problem. An important role is played by a Krein space realization result of
  Dijksma, Langer and de Snoo for functions analytic in a neighborhood of the origin. A key observation is that, using a simple Moebius transform,
  one can reduce the study of discrete analytic functions in the upper right quadrant to
  problems of function theory in the open unit disk. As an example, we associate to each finite positive measure on $[0,2\pi]$ a discrete analytic function on the right-upper
  quarter plane with values on the lattice defining a positive definite function. Emphasis is put on the rational case, both when an underlying Carath\'eodory function is rational and when, in the
  positive case, the spectral function is rational. The rational case and the general case are linked via the existence of a unitary dilation, possibly in  a Krein space.
  \end{abstract}
\tableofcontents


\date{today}
\tableofcontents
\section{Introduction}
\setcounter{equation}{0}
This introduction gives an overview of the paper and its background, and is divided into four parts.

\subsection{Prologue}
In this paper we study new aspects of discrete analytic functions, and connections with moment problems, structured matrices and reproducing kernel spaces of pairs. We propose a general scheme for
constructing  discrete analytic functions of the form
\begin{equation}
\label{newform123}
f(m,n)=CA_1^mA_2^nB,\quad m,n=0,1,\ldots
\end{equation}
where $A_1,A_2,B,C$ are operators between Hilbert spaces, or matrices of appropriate sizes, $A_1$ and $A_2$ satisfying moreover the relation
\begin{equation}
  \label{newconditioncom11}
I+iA_1-iA_2-A_1A_2=0,
\end{equation}
and their connections with discrete analytic functions developed in \cite{alpay2021discrete,alpay2021discrete2}.
Formulas \eqref{newform123}-\eqref{newconditioncom11} give the recipe  for a large number of explicit examples of discrete analytic functions.
Rewriting \eqref{newconditioncom11} as
\begin{equation}
  (A_1+iI)(A_2-iI)=2I,
    \end{equation}
we see that $A_1$ is completely determined by $A_2$ (and vice-versa) when these two operators commute. We will also be interested in expressions of the form

\begin{equation}
\label{newform12345}
f(m,n)=CA_1^mA_1^{*n}C^*,\quad m,n=0,1,\ldots
\end{equation}
This forces $A_1$ to solve the equation
\begin{equation}
  I+i(A_1-A_1^*)-A_1A_1^*=0.
  \label{s-a-+c-t-u}
\end{equation}
\begin{remark}{\rm
    Expressions of the type \eqref{s-a-+c-t-u} with right hand side equal to $0$, or ``close to $0$'', meaning that it has small rank,
    \begin{equation}
      I+i(A_1-A_1^*)-A_1A_1^*=CJC^*
    \end{equation}
    where $J\in\mathbb C^{m\times m}$ satisfies $J=J^2=J^{-1}$ and $C$ is a linear bounded operator from the underlying Hilbert space into
    $\mathbb C^m$,   form a recurring theme in this work;
see \eqref{equa-daf2} and \eqref{4-28-000}.}
\end{remark}

Another problem we will consider is to start from the values of the discrete analytic function on the positive axis, and assumed of the form
\begin{equation}
f(m,0)=CA^mC^*,\quad m=0,1,\ldots
\label{real-pb}
\end{equation}
To extend the function to a function of the form \eqref{newform12345} will make use the unitary dilation
theorems of Sz-Nagy and Foias \cite{nf} and of
Ch. Davis \cite{MR264438}.\smallskip

We note that realization theory of rational and analytic functions play an important role in the paper. Two different approaches and sets of tools are used. First, realizations where the
main operator is unitary. Such realization have an infinite dimensional state space in general.
Next, realizations of rational functions. Then, the state space is finite dimensional.
The two approaches are reconciled via the above mentioned unitary dilation theorems.

\begin{remark}{\rm In the paper we use a number of well-known methods and formulas from the theory of rational matrix-valued functions.
    We also use a number of classical definitions and results from operator theory (such as the existence of a unitary dilation for a contraction in Hilbert space) and from
    the theory of Pontryagin and Krein spaces.
    Having in view various audiences not necessarily familiar
        with these methods we have repeated a number of definitions and arguments. }
    \end{remark}

    \begin{remark} {\rm Recovering the pair $(C,A)$ in \eqref{real-pb} from all the values $f(m,0)$ is a special case of the realization problem in linear system theory;
        if only the values $f(0,0),\ldots, f(b,0)$ are given for some $b\in\mathbb N$, this is
        a particular case of the partial realization problem (see \cite{gkl-scl,MR669917}).}\end{remark}

\subsection{Discrete analytic functions and displacement structure}
Let  $\Lambda$ denote the lattice $\Lambda=\mathbb Z+i\mathbb Z$, and let $\Omega$ be a sub-lattice (i.e. $\Omega$ is a union of squares with vertices
$z,z+1,z+i,z+i+1$ with $z\in\mathbb Z$). A function $f(z)$ defined in $\Omega$ is called {\sl discrete analytic} if on every such square it holds that
\begin{equation}
  \label{dafdaf}
  \frac{f(z+i+1)-f(z)}{1+i}=\frac{f(z+1)-f(z+i)}{1-i}.
\end{equation}

Equivalently, and using the notation $f(m+in)=f(m,n)$, one can rewrite \eqref{dafdaf} as
\begin{equation}
  \label{def-daf}
  f(m,n)+if(m+1,n)-if(m,n+1)-f(m+1,n+1)=0,
\end{equation}
which is the discrete counterpart of the Cauchy-Riemann equations.\smallskip

These functions were introduced by Ferrand in \cite{MR0013411} and studied in
particular by Duffin in \cite{MR0078441}. Further studies include \cite{MR2195410,MR450576, MR0435417}.
We note that the above definitions still make sense for matrix-valued, or even operator-valued, functions. In view of \eqref{dafdaf}, finite sections
$(f(m,n))_{m,n=0}^{N}$ satisfy
\begin{equation}
  \delta f(m,n)  +\overline{\beta} f(m+1,n)+\beta f(m,n+1)+\alpha f(m+1,n+1)=0,\quad m,n=0,\ldots N,
  \label{heinig-rost-345}
\end{equation}
with
\[
\alpha=-1,\quad \beta=-i,\quad \delta=1,
  \]
and have a displacement structure in the sense of Kailath, Kung and Morf \cite{zbMATH03649904} and Heinig and Rost \cite{heinig-rost-laa}; see \cite[p. 147]{MR1197502} for a discussion of the latter.
\smallskip

As we explain in the sequel, equation \eqref{def-daf} allows to translate properties of the underlying discrete analytic function in terms of properties of associated reproducing kernel spaces related to the ones introduced by
de Branges and Rovnyak (see \cite{  dbr1, dbr2,MR0229011}). A first step in this direction is Zeilberger's formula for the generating function associated to the values of a function discrete analytic
in $\Lambda_{++}=\mathbb N_0+i \mathbb N_0$.
\subsection{Generating function and Zeilberger's formula}
Consider the generating function
\begin{equation}
  \label{richelieu-drouot}
  k_f({\lambda},{\nu})=\sum_{m,n=0}^\infty f(m,n){\lambda}^m\overline{{\nu}}^n
\end{equation}
of a $\mathbb C^{p\times q}$-valued discrete analytic function in $\Lambda_{++}$, which is invariant under the transformations
\begin{equation}
  \label{shifts}
  z\mapsto z+1\quad{\rm and}\quad z\mapsto z+i.
  \end{equation}
Sometimes we will identify $\Lambda_{++}$ and
$\mathbb N_0\times \mathbb N_0$.
Assuming  convergence in a neighborhood of the origin of $\mathbb C^2$, it follows from \cite[Theorem 1, p. 1247]{aro3} that there is a reproducing
kernel Hilbert space of pairs (see Definition \ref{groningen} below) with reproducing kernel \eqref{richelieu-drouot}. At this stage, the discrete analytic structure has not been taken into account.
To include the latter in the analysis we use Zeilberger's formula for the generating function
\begin{equation}
  \label{opera-bastille}
  k_f({\lambda},{\nu})=\frac{\Phi_L({\lambda})+\Phi_R({\nu})^*}{1+i{\lambda}-i\overline{{\nu}}-{\lambda}\overline{{\nu}}}
\end{equation}
with
\begin{eqnarray}
  \label{q1}
  \Phi_L({\lambda})&=&(1+i{\lambda})\left(\sum_{m=0}^\infty f(m,0){\lambda}^m\right)
               -\frac{f(0,0)}{2}+iX\\
  \label{q2}
  \Phi_R({\lambda})&=&(1+i{\lambda})\left(\sum_{n=0}^\infty f(0,n)^*{\lambda}^n\right)-\frac{f(0,0)^*}{2}+iX^*.
\end{eqnarray}
where $X\in\mathbb C^{p\times q}$ is an arbitrary matrix; see \cite[(3.20), p. 352]{MR450576}. Note that the latter paper considers more generally the case of formal power series.
The functions $\Phi_L$ and $\Phi_R$ fix the values of the function $f$ on the horizontal right axis and vertical upper axis respectively, which are in general independent one from the other.
  The formal power series given by the generating function is convergent in a neighborhood of the origin of $\mathbb C^2$ if and only if both $\Phi_L$ and $\Phi_R$ converge
  near the origin.\\

\begin{definition}
\label{left-right-gen}
The functions $\Phi_L$ and $\Phi_R$ will be called respectively  the left and right boundary generating functions.
\end{definition}

The denominator in \eqref{opera-bastille} can be rewritten as
\begin{equation}
1+i{\lambda}-i\overline{\nu}-{\lambda}\overline{{\nu}}=a({\lambda})\overline{a({\nu})}-b({\lambda})\overline{b({\nu})}
\end{equation}
with
\begin{equation}
  \label{opera-garnier}
    a({\lambda})=1+i{\lambda}\quad{\rm and}\quad b({\lambda})=\sqrt{2}{\lambda}.
  \end{equation}
  The kernel \eqref{opera-bastille}
can thus be rewritten as

\begin{equation}
  \label{opera-bastille111}
  k_f({\lambda},{\nu})=\frac{\Phi_L({\lambda})+\Phi_R({\nu})^*}{a({\lambda})\overline{a({\nu})}-b({\lambda})\overline{b({\nu})}},
\end{equation}
and  is therefore of the form of the kernels studied in a series of papers which includes \cite{ad-laa4,MR1197502,ad-jfa}. Such kernels are closely related to matrices
having a displacement structure and to a host of problems in operator theory and analytic functions.
Following these papers, we set
\begin{equation}
\label{sigma}
  \sigma({\lambda})=\frac{b({\lambda})}{a({\lambda})}=\frac{\sqrt{2}{\lambda}}{1+i{\lambda}}
\end{equation}
and
\begin{equation}
  \label{1-20-20-20}
    \begin{split}
      \Omega_+&=\left\{{\lambda}\in\mathbb C\,;\,|\sigma({\lambda})|<1\right\}\,=B(-i,\sqrt{2})\\
        \Omega_-&=\left\{{\lambda}\in\mathbb C\,;\,|\sigma({\lambda})|>1\right\}\,=E(-i,\sqrt{2})\\
                \Omega_0&=\left\{{\lambda}\in\mathbb C\,;\,|\sigma({\lambda})|=1\right\}\,=C(-i,\sqrt{2}),
                \end{split}
              \end{equation}
              where $E(-i,\sqrt{2})$ denotes the exterior of $B(-i,\sqrt{2})$. We can then rewrite \eqref{opera-bastille111} as

\begin{equation}
  \label{opera-bastille111111}
  k_f({\lambda},{\nu})=\frac{\Phi_L({\lambda})+\Phi_R({\nu})^*}{a({\lambda})(1-\sigma(\lambda)\overline{\sigma(\nu)})\overline{a({\nu})}}.
\end{equation}

               \begin{remark}{It is also important, especially when discussing connections with Schur analysis, to introduce the functions
                    \begin{equation}
                      \label{chphi}
                      \varphi_L(\lambda)=\Phi_L(\sigma^{-1}(\lambda))\quad{\rm and}\quad                      \varphi_R(\lambda)=\Phi_R(\sigma^{-1}(\lambda)).
\end{equation}
}
                \end{remark}

                \begin{definition}
                  \label{chara123}
                  We will call $\varphi_L$ and $\varphi_R$ given by
                  \eqref{chphi} the left and right characteristic
                  functions associated to $f$. When $\varphi_L=\varphi_R=\varphi$
                  we will call $\varphi$ the characteristic function associated to $f$.
                  \end{definition}

We note that the Moebius transformation $\sigma({\lambda})$ maps the open disk $B(-i,\sqrt{2})$ onto the open unit disk $B(0,1)$, with inverse $\sigma^{-1}({\lambda})=\frac{{\lambda}}{\sqrt{2}-i{\lambda}}$. The above mentioned papers considered  general pairs $(a,b)$ of analytic functions in a connected set $\Omega$, and for which the corresponding sets $\Omega_\pm$ (and hence $\Omega_0$)
are not empty.\smallskip

The case where  $f$ is symmetric (see Definition \ref{symdef}),
meaning that
\begin{equation}
  \label{symm123}
    f(m,n)=f(n,m)^*,\quad m,n=0,1,\ldots
\end{equation}
  is of special interest. Note that every discrete analytic function can be written as a linear combination of two
  discrete analytic symmetric functions as
\begin{equation}
  f(m,n)=\frac{f(m,n)+f(n,m)^*}{2}+i\frac{f(m,n)-f(n,m)^*}{2i}.
\end{equation}

For a symmetric discrete analytic function, we have $\Phi_L({\lambda})=\Phi_R({\lambda})\stackrel{\rm def.}{=}\Phi({\lambda})$. Since $\sigma^{-1}$ is analytic in a neighborhood of the origin the
function $\Phi(\sigma^{-1}({\lambda}))$ has the same property, and by a result of Dijksma, Langer and de Snoo (see \cite[Theorem 1, p. 126]{MR903068}; Theorem \ref{dlds} below), we can write the characteristic function in the form
\begin{equation}
  \label{gare-de-l-est}
  \Phi(\sigma^{-1}({\lambda}))=iX+\frac{1}{2}C(U+{\lambda}I)(U-{\lambda}I)^{-1}C^*,
  \end{equation}
  where $U$ is a  bounded unitary operator in a Krein space $\mathfrak K$ and $C$ is a bounded operator from $\mathfrak K$ into $\mathbb C^n$ (see Definition \ref{def-krein} for
  the notion of Krein space; the results of \cite{MR903068} are much more general and  are given in the framework of operator valued functions; we focus here on the matrix-valued case). In \eqref{gare-de-l-est}, $C^*$ denotes the adjoint of $C$ with respect to the indefinite inner product in $\mathfrak K$.
  The corresponding discrete analytic function is given by
\begin{equation}
      \label{gare-du-nord1}
      f(m,n)=C(\sqrt{2}U^{-1}-iI)^m(\sqrt{2}U+iI)^nC^*,\quad n,m=0,1,\ldots
    \end{equation}

    One of the advantages of the above formula is that it allows to characterize discrete analytic extensions in terms of the spectrum of $U$.
 A unitary operator in a Pontryagin space, and more generally in a Krein space, can have spectrum outside the unit circle. When the points $-i\sqrt{2}$ and $\frac{-i}{\sqrt{2}}$
      are not in the
      spectrum of $U$, the function $f$ has a unique discrete analytic extension to $\Lambda$, given by \eqref{gare-du-nord1}. \smallskip

    We furthermore note that \eqref{symm123}
implies the existence of complex numbers $c_{m,u}$ such that
 \[
f(m,n)=\sum_{a}c_{m,a}M_a+\sum_{b} \overline{c_{n,b}}M_b^*,
\]
where $M_b=CU^bC^*$ denote the ``moments'' of the pair $(C,U)$.  This remark is of particular importance in the Hilbert space setting, where the $M_b$ are moments of an underlying measure on
$[0,2\pi]$; see Theorem \ref{th1} below.\\

 Two cases of special interest of \eqref{gare-du-nord1} are when $\mathfrak K$ is a Pontryagin space, or a Hilbert space. In the first case, the finite
 section matrices have uniformly at most a given  finite number of negative squares. In the Hilbert space case they are positive. Equivalently, but still under the hypothesis of convergence of
 the generating function in a neighborhood of the origin, the kernel
 \begin{equation}
   L_\Phi(\lambda,\nu)=\frac{\Phi({\lambda})+\Phi(\nu)^*}{1+i{\lambda}-i\overline{\nu}-{\lambda}\overline{\nu}}
      \end{equation}
      has a finite number of negative squares in the first case, and is positive definite in the second case (see
      \cite[Lemma 1.1.6, p. 10]{adrs} for the positive case).
      In that case, and as explained after Corollary \ref{mumomu} below,  the operator $U$ in \eqref{gare-de-l-est} and \eqref{gare-du-nord1} can be chosen to be multiplication by $e^{it}$ in
      an associated Lebesgue space of the unit circle. The function \eqref{gare-de-l-est} is then respectively a generalized Carath\'eodory function and a Carath\'eodory function;
      see Definition \ref{cara789} for the latter.   \smallskip


\subsection{Discrete analytic functions and moments problems}
The case $\delta=-\alpha=1$ and $\beta=0$ in \eqref{heinig-rost-345} corresponds to Toeplitz
matrices and block Toeplitz matrices. Positive block Toeplitz matrices (say, with blocks in $\mathbb C^{p\times p}$)
are characterized in terms of
moments of a $\mathbb C^{p\times p}$-valued positive measure on $[0,2\pi]$. More precisely (see for instance \cite{akhiezer,knud,
MR0008438} for the scalar case and \cite{ggk1,ggk2} and references therein for the matrix-valued case):

\begin{theorem}
  Let $T_0,T_1,\ldots\in\mathbb C^{p\times p}$. There exists a finite $\mathbb C^{p\times p}$-valued positive measure
  $M$ on $[0,2\pi]$ such that
  \begin{equation}
T_j=\frac{1}{2\pi}\int_0^{2\pi}e^{-ijt}dM(t),\quad j=0,1,\ldots
  \end{equation}
  if only if all the block Toeplitz matrices
  \begin{equation}
    \mathbf T_N=\begin{pmatrix} T_0&T_{1}^*&T_2^*&\cdots & T_{N}^*\\
    T_1&T_0&T_1^*&\cdots &T_{N-1}^*\\
    \ddots&\ddots& \ddots& &\ddots\\
    T_N&T_{N-1}&\cdots &\cdots&T_0
    \end{pmatrix},\quad N=0,1,\ldots
    \end{equation}
are non-negative.
\end{theorem}
In a similar way we prove the following, which is one of the main results in the paper.

\begin{theorem}
  \label{th1}
  Let $f(z)$ be a $\mathbb C^{p\times p}$-valued discrete analytic symmetric function in the quarter plane $\Lambda_{++}$.
  All the finite section matrices $F_N=(f(m,n))_{m,n=0}^{N}$ are non-negative if and only if there exists a
  $\mathbb C^{p\times p}$-valued finite positive measure $M(t)$ on $[0,2\pi]$ such that
    \begin{equation}
\label{formulemoments}
      f(m,n)=\frac{1}{2\pi}\int_0^{2\pi}(\sqrt{2}e^{-it}-i)^m(\sqrt{2}e^{it}+i)^ndM(t),\quad m,n=0,1\ldots
  \end{equation}
\end{theorem}

Because of connections with operator theory, linear systems and inverse problems, we will call $M$ the spectral measure associated to the function $f(m,n)$. Formula \ref{formulemoments}
suggests connections with the work \cite{ab-1994}. These will be considered elsewhere.\smallskip

As a consequence of \eqref{formulemoments} we have:

\begin{corollary}
  \label{mumomu}
Under the hypothesis of the previous theorem, the function $f$ has a discrete analytic extension to $\Lambda$, still given by \eqref{formulemoments}. This extension is unique.
  \end{corollary}

Equation \eqref{formulemoments} is a special case of \eqref{gare-du-nord1},
  where $Uf(t)=e^{it}f(t)$ is the unitary operator of multiplication by $e^{it}$ in $\mathbf L^2([0,2\pi],dM)$, and $Cx=\int_o^{2\pi}dM(t)x(t)$ is linear and bounded from $\mathbf L^2([0,2\pi],dM)$
  into $\mathbb C^p$.

  \begin{definition}
    \label{dafm}
    We call a matrix $(f(m,n))_{n,m=0}^N$ whose entries are the values of a discrete analytic functions at the point $(m,n)$, $m,n=0,\ldots, N$ a matrix
    with discrete analytic structure. It will be called a block matrix with discrete analytic structure if the function is matrix-valued.
  \end{definition}

  When the matrices are non-negative, the correspondence between matrices with discrete analytic structure and Toeplitz matrices allows to consider various extension problems
  (such as the one-step extension procedure and maximum entropy analysis; see e.g. \cite{van-den-bos}) in the present setting.

  \smallskip

  The rational case is of special interest, and is considered in two cases. When the function
  $\Phi$ above is rational, or when the spectral measure $M(t)$ in \eqref{formulemoments} is absolutely
  continuous with respect to Lebesgue measure, with rational density function. In both cases,
  we use realization theory to get formulas for $f(m,0)$. We obtain formulas of the kind
  \[
  f(m,0)=\mathsf C(\sqrt{2}\mathsf A-iI_N)^m\mathsf D
  \]
  where $\mathsf A\in\mathbb C^{N\times N}$ and $\mathsf C$ and $\mathsf D$ are matrices of compatible sizes. $\mathsf A$ is in general not unitary, whether in a Hilbert or Krein space metric, and
  one needs to use a unitary dilation of $A$ to get a formula of the form \eqref{newform123} for the extension of $f(m,0)$ to $f(m,n)$.\smallskip

  In this paper we consider matrix-valued functions, but the operator-valued case is also
  of interest. Herglotz integral formula holds in the Hilbert space setting
    (see \cite[Theorem 4.5, p. 23]{MR48:904}) and in the setting of functions taking values from a topological vector space into its anti-dual; see \cite{atv1,atv2}. One can therefore extend the
    results of the present paper to these settings, which play an important role in stochastic processes (see e.g. \cite{lifs,MR626346}).\\

    \subsection{Outline of the paper} The paper consists of four sections, besides this Introduction, each of which is divided into subsections to improve the
    readability, as described in the table of contents. Section 2 contains some preliminaries on rational functions, also in the discrete case, and their realizations. We also discuss Krein and Pontryagin spaces, and structured matrices which are relevant  for reproducing kernel pairs of spaces in duality. In Section 3 we prove various results for realizations of the the coefficients symmetric discrete analytic functions in the Krein and Pontryagin spaces. The positive case is considered in Section 4. Here, formulas for the coefficients of a symmetric discrete analytic function are given using Herglotz integral formula. A link with moment problems is also presented. Finally, in Section 5 we consider the rational case.

\section{Preliminaries}
  \setcounter{equation}{0}
\subsection{Realization of rational functions}
\label{real-sec}
We recall that any matrix-valued rational function $R({\lambda})$ analytic at infinity can be written in the form
\begin{equation}
  \label{real}
  R({\lambda})=D+C({\lambda}I_N-A)^{-1}B
\end{equation}
Expression \eqref{real} is called a realization of $R({\lambda})$. It is minimal when the size $N$ of the matrix $A$ is minimal, and then is unique up to a similarity matrix, meaning that if
$R({\lambda})=D_j+C_j({\lambda}I_N-A_j)^{-1}B_j$, $j=1,2$, are two minimal realizations of the $\mathbb C^{p\times q}$-valued rational function $R$, then $D_1=D_2=D=R(\infty)$ and there exists an invertible similarity matrix $T$,
uniquely determined by the two realizations, and such that
\[
\begin{pmatrix}T&0\\0&I_p\end{pmatrix}\begin{pmatrix}A_1&B_1\\ C_1&D\end{pmatrix}=\begin{pmatrix}A_2&B_2\\ C_2&D\end{pmatrix}\begin{pmatrix}T&0\\0&I_q\end{pmatrix}.
\]

Minimality of the realization \eqref{real} is equivalent to the pair $(C,A)$  being observable and the pair $(A,B)$ being controllable, meaning respectively
\begin{equation}
  \label{obs}
  \bigcap_{n=0}^\infty \ker CA^n=\left\{0\right\}
\end{equation}
and
\begin{equation}
  \label{contro}
  \bigcup_{n=0}^\infty {\rm ran}\, A^nB=\mathbb C^N.
\end{equation}
See for instance \cite{bgk1,MR0255260}. One then says that the triple $(A,B,C)$ is minimal.\smallskip

As a first illustration of the above notions we present following proposition, to be used in Remark \ref{followlem}.

 \begin{proposition}
    \label{lemma31}
    Assume the triple $(A,B,C)$ minimal and both $A$ and $I_N+iA$ invertible. Then the triple
\begin{equation}
    \label{trip}
  (A^{-1}, C(I_N+iA)^{-1}, A^{-1}B)
\end{equation}
  is also minimal.
\end{proposition}

\begin{proof}
We first prove that the pair $(C(I_N+iA)^{-1}, A^{-1})$ is observable. For $t\in \mathbb C$ small enough  we can write
\begin{equation}
  \label{sumsum}
\sum_{a=0}^\infty t^aC(I_N+iA)^{-1}A^{-a}=CA(A-tI_N)^{-1}(I_N+iA)^{-1}.
\end{equation}
Let $\xi\in\cap_{a=0}^\infty \ker C(I_N+iA)^{-1}A^{-a}$. Then, \eqref{sumsum} implies that $CA(A-tI_N)^{-1}(I_N+iA)^{-1}\xi\equiv 0$ near the origin. Since the pair $(C,A)$ is observable, this implies that
$(I_N+iA)^{-1}\xi=0$ and hence $\xi=0$.\smallskip

One proves in a similar way that  the pair $(B^*A^{-*}, A^{-*})$ is observable. Hence the pair
$(A^{-1},A^{-1}B)$ is controllable, and the triple \eqref{trip} is minimal.
  \end{proof}

  We now consider realizations in the symmetric case, $\Phi_L=\Phi_R^*$. {\sl A priori} one has
  $f(m,0)=C_1A_1^mB_1$ and $f(0,n)=C_2A_2^nB_2$ for different triples $(A_1,B_1,C_1)$ and
  $(A_2,B_2,C_2)$ and one has
  \[
C_1A_1^mB_1=B_2^*A_2^{*m}C_2,\quad m=0,1,\ldots
\]
When the two triples are minimal there is uniquely determined similarity matrix $S$ such that
\[
C_1=B_2^*S,\quad A_1=S^{-1}A_2^{*}S,\quad B_1=S^{-1}C_2^*
\]
and little more seems to be derived from this; one cannot deduce from the uniqueness of the similarity matrix that it is (for instance) Hermitian. On the other hand, one has:
\begin{proposition}
Assume in the above $C_1=C_2=C$ and $B_1=B_2=B$. Then $S$ is Hermitian, and $A_1$ is Hermitian with respect to the indefinite metric defined by $S^{-1}$ (see \eqref{adj} for the latter).
\end{proposition}

\begin{proof}
By minimality  we have
  \[
    B^*=CS,\quad A_1^*=S^{-1}A_1S,\quad C^*=S^{-1}B
  \]
  These equations are also satisfied by $S^*$ and so $S=S^*$. The condition $A_1^*=S^{-1}A_1S$ then means that $A_1$ is self-adjoint with respect to the (possibly) indefinite metric defined by $S$.
  \end{proof}

We will need (see the proof of Proposition \ref{amsterdam}) the formula for the realization of a  product of two matrix-valued
rational functions analytic at infinity and of compatible sizes in terms of their respective realizations. More precisely (see e.g. \cite{bgk1,MR0325201}):
\begin{equation}
  \label{formprod}
(D_1+C_1({\lambda}I_{N_1}-A_1)^{-1}B_1)(D_2+C_2({\lambda}I_{N_2}-A_2)^{-1}B_2)=D+C({\lambda}I_N-A)^{-1}B
  \end{equation}
with $D=D_1D_2$ and
\begin{equation}
A=\begin{pmatrix}A_1&B_1C_2\\0&A_2\end{pmatrix},\quad B=\begin{pmatrix}B_1D_2\\ B_2\end{pmatrix}\quad {\rm and}\quad C=\begin{pmatrix}C_1&D_1C_2\end{pmatrix}.
\end{equation}

Realizations of the kind \eqref{real} are called ``centered at infinity''. In the sequel realizations centered at the origin, of the form
\begin{equation}
  R({\lambda})=D+{\lambda}C(I_N-{\lambda}A)^{-1}B
  \label{real-0}
\end{equation}
also appear. We note the power series expansion
\begin{equation}
  R({\lambda})=D+\sum_{u=1}^\infty {\lambda}^uCA^{u-1}B
  \label{exanR}
\end{equation}
in a neighborhood of the origin for such a realization.\smallskip

Every rational function with no pole at the origin can be written in the form \eqref{real-0}. Starting from a realization centered at infinity and assuming $A$ invertible in \eqref{real} we have
\begin{equation}
  \label{real-zero}
  \begin{split}
    R({\lambda})&=D+C({\lambda}I_N-A)^{-1}B \\
    &=D-CA^{-1}B+C\left\{({\lambda}I_N-A)^{-1}+A^{-1}\right\}B\\
    &=D-CA^{-1}B-{\lambda}CA^{-1}(I_N-{\lambda}A^{-1})^{-1}A^{-1}B
  \end{split}
  \end{equation}
  making the connections between realizations centered at $0$ and $\infty$. Other kind of realizations are possible, centered at an arbitrary point different from the origin,
  but will not be considered here. See \cite{gk1}.\smallskip

For the next theorem, see \cite[Theorem 4.7, p. 215]{ag}, \cite{MR2663312} for the self-adjoint case, one has to multiply by
$-i$ to get to the present case of functions taking skew self-adjoint values on the unit circle. A proof is provided for completeness.


  \begin{theorem}
    Let $\varphi$ be a $\mathbb C^{p\times p}$-valued rational function, analytic at infinity, and let $\varphi({\lambda})=D+C({\lambda}I_N-A)^{-1}B$ be a minimal realization of $\varphi$. Then
    $\varphi$ is skew self-adjoint on the unit circle, i.e.
    \[
\varphi({\lambda})+\varphi({\lambda})^*=0,\quad |{\lambda}|=1
\]
if and only if:\\
$(1)$ $\varphi$ is analytic at the origin (and so $A$ is invertible).\\
$(2)$ There exists a $\mathbb C^{N\times N}$ Hermitian invertible matrix $H$ such that

\begin{align}
  \label{a-unit}
  A&=H^{-1}(A^*)^{-1}H\\
  \label{a-unit-2}
  D+D^*&=-CH^{-1}C^*\\
  B&=-AH^{-1}C^*.
     \label{a-unit3}
\end{align}
Then, it holds that
\begin{equation}
  \label{comp1}
\frac{\varphi({\lambda})+\varphi({\nu})^*}{1-{\lambda}\overline{{\nu}}}=C({\lambda}I_N-A)^{-1}H^{-1}({\nu}I_N-A)^{-*}C^*,\quad {\lambda},{\nu}\in\rho(A),
\end{equation}
where $\rho(A)$ denotes the resolvent set of $A$, and
\begin{equation}
  \label{comp22}
  \varphi({\lambda})=\frac{D-D^*}{2}+\frac{1}{2}C(A+{\lambda}I_N)(A-{\lambda}I_N)^{-1}H^{-1}C^*.
\end{equation}
\label{th-ag-p-215}
\end{theorem}

\begin{proof}
  The necessity is based on the uniqueness of the minimal realization up to a similarity matrix. We have
  \begin{equation}
    \varphi(\lambda)=-(\varphi(1/\overline{\lambda}))^*
    \label{equa-phi}
    \end{equation}
  when defined. From this equation we deduce that analyticity at infinity implies analyticity at the origin, and in particular $A$ is invertible. Using \eqref{real-zero} we then
  rewrite \eqref{equa-phi}
  as
  \begin{equation}
    \label{realstar123}
D+C(\lambda I_N-A)^{-1}B=-D^*+B^*A^{-*}C^*+B^*A^{-*}(\lambda I_N-A^{-*})^{-1}A^{-*}C^*,
\end{equation}
from which we get
\begin{equation}
  \label{R-infty}
  D=-D^*+B^*A^{-*}C^*.
  \end{equation}
  \eqref{realstar123} is an equality between two minimal realizations of a given rational function, and hence there is a unique invertible similarity matrix, denoted here by $-H^{-1}$ such that
  \begin{align}
    \label{p111}
    C&=-B^*A^{-*}H,\\
    \label{p555}
    A&=H^{-1}A^{-*}H,\\
    B&=-H^{-1}A^{-*}C^*.
       \label{p222}
    \end{align}
    These equations are also satisfied by $-H^{-*}$, and by uniqueness of the similarity matrix we get $H=H^*$. Equations \eqref{a-unit}-\eqref{a-unit3} follow then from \eqref{R-infty} and
    \eqref{p111}-\eqref{p222}.\smallskip

We now prove \eqref{comp1} and \eqref{comp22} when \eqref{a-unit}-\eqref{a-unit3} are in force.
 We have
  \[
  \begin{split}
    \varphi({\lambda})+\varphi({\nu})^*&=D+D^*+C({\lambda}I_N-A)^{-1}B+B^*(\overline{{\nu}}I_N-A^*)^{-1}C^*\\
    &=-C{H^{-1}}C^*-C({\lambda}I_N-A)^{-1}A{H^{-1}}C^*-C{H^{-1}}A^*(\overline{{\nu}}I_N-A^*)^{-1}C^*\\
    &=C({\lambda}I_N-A)^{-1}\left\{-({\lambda}I_N-A){H^{-1}}(\overline{{\nu}}I_N-A^*)-\right.\\
    &\hspace{5mm}\left.-A{H^{-1}}(\overline{{\nu}}I_N-A^*)-({\lambda}I_N-A){H^{-1}}A^*\right\}(\overline{{\nu}}I_N-A^*)^{-1}C^*\\
    &=C({\lambda}I_N-A)^{-1}\left\{-{\lambda}\overline{{\nu}}{H^{-1}}+A{H^{-1}}A^*\right\}(\overline{{\nu}}I_N-A^*)^{-1}C^*\\
    &=(1-{\lambda}\overline{{\nu}})C({\lambda}I_N-A)^{-1}(\overline{{\nu}}I_N-A^*)^{-1}C^*,
  \end{split}
\]
where in particular we used \eqref{a-unit} (that is,
$A{H^{-1}}A^*={H^{-1}}$) to go from the second-to-last line to the last line.\smallskip

To prove \eqref{comp22} we write
\[
  \begin{split}
    \varphi({\lambda})&=\frac{D-D^*}{2}+\frac{D+D^*}{2}+C({\lambda}I_N-A)^{-1}B\\
    &=\frac{D-D^*}{2}-\frac{1}{2}C{H^{-1}}C^*-C({\lambda}I_N-A)^{-1}A{H^{-1}}C^*\\
        &=\frac{D-D^*}{2}-\frac{1}{2}C\left\{({\lambda}I_N-A)+2A\right\}({\lambda}I_N-A)^{-1}{H^{-1}}C^*,
    \end{split}
  \]
from which the formula follows.
  \end{proof}
We note that \eqref{a-unit} means that $A$ is unitary with respect to the indefinite metric induced by $H$ in $\mathbb C^N$:
\begin{equation}
  \label{metric-H}
[x,y]_H=y^*Hx,\quad x,y\in\mathbb C^N.
\end{equation}
Formula \eqref{comp1} expresses the fact that $H$ is the Gram matrix of the space spanned
by the columns of the function $C(\lambda I_N-A)^{-1}$.\smallskip

We refer to \cite{MR2186302} for linear algebra in presence of an indefinite metric. Denoting by $\mbox{}^{[*]}$ the adjoint with respect to the metric defined by $H$, i.e.
\begin{equation}
  \label{adj}
  A^{[*]}=H^{-1}AH,
  \end{equation}
we have from \eqref{p555} that $AA^{[*]}=I_N$, that is, the operator $A$ is unitary in the metric \eqref{metric-H}:
\[
[Ax,Ay]_H=[x,y]_H,\quad x,y\in\mathbb C^N.
\]

Re-expressing \eqref{comp22} and \eqref{comp1} in terms if the metric \eqref{metric-H} we have:

\begin{corollary}
  Formula \eqref{comp22} can be rewritten as
  \begin{equation}
    \varphi({\lambda})=\frac{D-D^*}{2}+\frac{1}{2}C(A+{\lambda}I_N)(A-{\lambda}I_N)^{-1}C^{[*]}.
  \end{equation}
  and \eqref{comp1} becomes
  \begin{equation}
    \label{compare-KYP}
    \frac{\varphi({\lambda})+\varphi({\nu})^*}{1-{\lambda}\overline{\nu}}=C({\lambda}I_N-A)^{-1}(\nu I_N-A)^{-[*]}C^{[*]},\quad \lambda,\nu\in\rho(A).
    \end{equation}
\end{corollary}

\begin{remark}{\rm
    When ${\rm Re}\,\varphi({\lambda})\ge 0$ in the open unit disk then (and only then) $H>0$ and can be chosen to be equal to $I_N$.
    Then, $\varphi$ is called a Carath\'eodory function. More generally, Carath\'eodory functions are functions analytic in the open unit disk and with a positive real part there
    See Definition \ref{cara789} below. Another important function is $S({\lambda})=(\varphi({\lambda})-I_p)(\varphi({\lambda})+I_p)^{-1}$, studied in the
    following two propositions.}
  \end{remark}

Formula \eqref{below} was obtained by one of the authors (D.A.) with Dr. I. Paiva.
See Acknowledgments section.

\begin{proposition}
\label{ap-prop}
Assume $D=D^*$ and $I_p=\frac{CH^{-1}C^*}{2}$ (that is, $\varphi(0)=I_p$).
Then
\begin{equation}
  \label{below}
(\varphi({\lambda})-I_p)(\varphi({\lambda})+I_p)^{-1}=\frac{{\lambda}}{2}CA^{-1}\left(I_N-{\lambda}\left(I_N-\frac{1}{2}HC^*C\right)A^{-1}\right)^{-1}H^{-1}C^*,
\end{equation}
or, in terms of the metric defined by $H$,
\begin{equation}
(\varphi({\lambda})-I_p)(\varphi({\lambda})+I_p)^{-1}=\frac{{\lambda}}{2}CA^{[*]}\left(I_N-{\lambda}\left(I_N-\frac{1}{2}C^{[*]}C\right)A^{[*]}\right)^{-1}C^{[*]}.
\end{equation}
\end{proposition}

\begin{proof} We can write
  \[
    \begin{split}
      (\varphi({\lambda})-I_p)(\varphi({\lambda})+I_p)^{-1}&=I_p-2(\varphi({\lambda})+I_p)^{-1}\\
      &=I_p-2(I_p+\frac{1}{2}C{H^{-1}}C^*+\frac{1}{2}C(A+{\lambda}I_N)(A-{\lambda}I_N)^{-1}{H^{-1}}C^*-\frac{1}{2}C{H^{-1}}C^*)^{-1}\\
      &=I_p-2(2I_p+\frac{1}{2}C\left\{(A+{\lambda}I_N)(A-{\lambda}I_N)^{-1}-I_N\right\}{H^{-1}}C^*)^{-1}\\
      &=I_p-2(2I_p+{\lambda}C(A-{\lambda}I_N)^{-1}{H^{-1}}C^*)^{-1}\\
      &=I_p-(I_p+\frac{1}{2}{\lambda}C(A-{\lambda}I_N)^{-1}{H^{-1}}C^*)^{-1}\\
      &=I_p-(I_p-\frac{1}{2}{\lambda}C(I_N+{\lambda}\frac{1}{2}(A-{\lambda}I_N)^{-1}{H^{-1}}C^*C)^{-1}(A-{\lambda}I_N)^{-1}{H^{-1}}C^*
      \end{split}
    \]
    where we have used $(I+MU)^{-1}=I-M(I+UM)^{-1}U$ for matrices $M$ and $U$ of appropriate sizes with $M=\frac{1}{2}{\lambda}C$ and $U=(A-{\lambda}I_N)^{-1}{H^{-1}}C^*$.
    Hence,
    \[
 \begin{split}
   (\varphi({\lambda})-I_p)(\varphi({\lambda})+I_p)^{-1}&=\frac{1}{2}{\lambda}C(I_N+{\lambda}\frac{1}{2}(A-{\lambda}I_N)^{-1}{H^{-1}}C^*C)^{-1}(A-{\lambda}I_N)^{-1}{H^{-1}}C^*\\
   &=\frac{{\lambda}}{2}C(A-{\lambda}I_N+\frac{{\lambda}}{2}{H^{-1}}C^*C)^{-1}{H^{-1}}C^*\\
   &=\frac{{\lambda}}{2}CA^{-1}\left(I_N-{\lambda}\left(I_N-\frac{1}{2}{H^{-1}}C^*C\right)A^{-1}\right)^{-1}{H^{-1}}C^*.
    \end{split}
    \]
  \end{proof}

  The function $S({\lambda})=(\varphi({\lambda})-I_p)(\varphi({\lambda})+I_p)^{-1}$ satisfies $S({\lambda})S({\lambda})^*=I_p$ on
  the unit circle. Thus the following result is a consequence of \cite[Theorem 3.1, p. 197]{ag}. We give a direct proof
  for completeness and verification.

\begin{proposition}
The matrix
  \[
    \begin{pmatrix}(I_N-\frac{1}{2}{{H^{-1}}C^*C})A^{-1}&\frac{1}{\sqrt{2}}{H^{-1}}C^*
      \\ \frac{CA^{-1}}{\sqrt{2}}&0\end{pmatrix}
\]
satisfies
\[
  \begin{split}
  \begin{pmatrix}(I_N-\frac{1}{2}{{H^{-1}}C^*C})A^{-1}&\frac{1}{\sqrt{2}}{H^{-1}}C^*
  \\ \frac{CA^{-1}}{\sqrt{2}}&0\end{pmatrix}
\begin{pmatrix}{H^{-1}}&0\\0&I_p\end{pmatrix}
\begin{pmatrix}(I_N-\frac{1}{2}{{H^{-1}}C^*C})A^{-1}&\frac{1}{\sqrt{2}}{H^{-1}}C^*
  \\ \frac{CA^{-1}}{\sqrt{2}}&0\end{pmatrix}^*&=\\
&\hspace{-2cm}=
\begin{pmatrix}{H^{-1}}&0\\0&I_p\end{pmatrix}.
\end{split}
  \]
\end{proposition}

\begin{proof}
Equation \eqref{a-unit} can be rewritten as $A^{-1}{H^{-1}}={H^{-1}}A^*$ and hence
  \[
    \begin{split}
      \frac{1}{\sqrt{2}}CA^{-1}{H^{-1}}A^{-*}(I_N-\frac{1}{2}C^*C{H^{-1}})&=-\frac{1}{\sqrt{2}}C{H^{-1}}(I_N-\frac{1}{2}C^*C{H^{-1}})\\
      &=\frac{1}{\sqrt{2}}C{H^{-1}}+
      \underbrace{\frac{1}{2}C{H^{-1}}C^*}_{=I_p}\frac{1}{\sqrt{2}}C{H^{-1}}\\
      &=0.
    \end{split}
  \]
  Still using \eqref{a-unit}, we have:
  \[
    \begin{split}
      (I_N-\frac{1}{2}{{H^{-1}}C^*C})A^{-1}{H^{-1}}A^{-*}(I_N-\frac{1}{2}C^*C{H^{-1}})+\frac{1}{2}{H^{-1}}C^*C{H^{-1}}&=\\
      &\hspace{-10cm}=
      (I_N-\frac{1}{2}{{H^{-1}}C^*C}){H^{-1}}(I_N-\frac{1}{2}C^*C{H^{-1}})+\frac{1}{2}{H^{-1}}C^*C{H^{-1}}\\
      &\hspace{-10cm}={H^{-1}}-\frac{1}{2}{H^{-1}}C^*C{H^{-1}}-\frac{1}{2}{H^{-1}}C^*C{H^{-1}}+\frac{1}{4}{H^{-1}}C^*C{H^{-1}}C^*C{H^{-1}}+\frac{1}{2}{H^{-1}}C^*C{H^{-1}}\\
      &\hspace{-10cm}={H^{-1}}
    \end{split}
  \]
  since $I_p=\frac{1}{2}C{H^{-1}}C^*$.
          \end{proof}
We note that a far-reaching generalization of \eqref{comp22} is given by the realization result of Dijksma-Langer-de Snoo stated below (see Theorem \ref{kl-extension} and formula \eqref{realcara}),
where $\varphi$ is an arbitrary operator-valued function analytic in a neighborhood of the origin, and the finite dimensional Pontryagin space $(\mathbb C^N,[\cdot,\cdot]_H)$ is replaced by a
Krein space (Pontryagin and Krein spaces are defined in Section \ref{krpon} below).\smallskip

  When the function $\varphi$ is rational and analytic at infinity, but not necessarily lossless, there is another far-reaching counterpart of
\eqref{comp22}, namely the Kalman-Yakubovich-Popov lemma (KYP lemma); see \cite{Anderson_Moore,faurre,MR525380} and \cite{DDGK}.
The latter gives a characterization of the minimal characterizations of $\varphi$. This realization
  in turn allows to get information on the coefficients of the underlying discrete analytic function.
  \begin{equation}
    \varphi({\lambda})=D+C({\lambda}I_N-A)^{-1}B
  \end{equation}
  and there exist matrices $L\in\mathbb C^{N\times r}$,
  $W\in\mathbb C^{p\times r}$ (where $r\in\mathbb N_0$)
  and a self-adjoint matrix $P\in\mathbb C^{N\times N}$ such that
  \begin{eqnarray}
    \label{kyp1}
    P- APA^*&=&-LL^*\\
    \label{kyp2}
    APC^*&=&-B-LW^*\\
    WW^*&=&D+D^*+CPC^*
            \label{kyp3}
  \end{eqnarray}

  \begin{remark} {\rm The difference of signs between the equations \eqref{kyp1}-\eqref{kyp3} stems from the fact that functions with a positive
        real part outside the closed unit disk are considered, while we consider here functions analytic in the open unit disk.}
          \end{remark}

   It is of interest to compare these equations to \eqref{a-unit}-\eqref{a-unit3}.
    We first remark that \eqref{a-unit} may be rewritten as $A^*(HAH^{-1})=I_N$, hence, $H(AH^{-1}A^*)=I_N$ and so
    \begin{equation}
      H^{-1}=AH^{-1}A^*.
      \label{a-unit!}
    \end{equation}
With $L=0$ and $W=0$ in \eqref{kyp1}-\eqref{kyp3} (i.e. $r=0$) we recover \eqref{a-unit-2}-\eqref{a-unit3} and \eqref{a-unit!} with $H^{-1}$ replaced by $P$.
The latter will be invertible when the original realization is minimal, making the connection with Theorem \ref{th-ag-p-215}.
In general $P$ is not unique in general and need not be invertible. The original KYP lemma considers the case where $P\ge0$; see the above mentioned references and in particular \cite{MR525380}. The case where general Hermitian matrices $P$ are allowed is
studied in \cite{DDGK}.\smallskip

 In the following proposition the formula \eqref{decomp123}  gives a proof of the sufficiency in the KYP lemma, and should be compared with \eqref{compare-KYP}.

  \begin{proposition}
    It holds that
    \begin{equation}
      \label{decomp123}
      \begin{split}
        \frac{\varphi({\lambda})+\varphi({\nu})^*}{1-{\lambda}\overline{\nu}}&=C({\lambda}I_N-A)^{-1}P(\overline{\nu}I_N-A^*)^{-1}C^*+\\
        &\hspace{5mm}+\frac{(W-C({\lambda}I_N-A)^{-1}L)(W^*-L^*(\overline{\nu}I_N-A^*)^{-1}C^*)}{1-{\lambda}\overline{\nu}},\quad \lambda,\nu\in\rho(A).
        \end{split}
\end{equation}
\end{proposition}
  \begin{proof}
Let us compute
\[
  \begin{split}
    \varphi({\lambda})+\varphi({\nu})^*-(1-{\lambda}\overline{\nu})C({\lambda}I_N-A)^{-1}P(\overline{\nu}I_N-A^*)^{-1}C^*&=\\
&\hspace{-80mm}   = WW^*-CPC^*-C({\lambda}I_N-A)^{-1}(LW^*+APC^*)-(WL^*+CPA^*)(\overline{\nu}I_N-A^*)^{-1}C^*-\\
  &\hspace{-75mm}-(1-{\lambda}\overline{\nu})C({\lambda}I_N-A)^{-1}P(\overline{\nu}I_N-A^*)^{-1}C^*\\
  &\hspace{-80mm}=WW^*-C({\lambda}I_N-A)^{-1}LW^*-WL^*(\overline{\nu}I_N-A^*)^{-1}C^*+\\
  &\hspace{-75mm}+C({\lambda}I_N-A)^{-1}\left\{\Delta\right\}(\overline{\nu}I_N-A^*)^{-1}C^*,
\end{split}
\]
with
\[
  \begin{split}
    \Delta&=-({\lambda}I_N-A)P(\overline{\nu}I_N-A^*)-AP(\overline{\nu}I_N-A^*)-({\lambda}I_N-A)PA^*-(1-{\lambda}\overline{\nu})P\\
    &=APA^*-P\\
    &=LL^*,
    \end{split}
  \]
  and so
  \[
    \begin{split}
WW^*-C({\lambda}I_N-A)^{-1}LW^*-WL^*(\overline{\nu}I_N-A^*)^{-1}C^*+
+C({\lambda}I_N-A)^{-1}\left\{\Delta\right\}(\overline{\nu}I_N-A^*)^{-1}C^*=\\
&\hspace{-155mm}=WW^*-C({\lambda}I_N-A)^{-1}LW^*-WL^*(\overline{\nu}I_N-A^*)^{-1}C^*+
C({\lambda}I_N-A)^{-1}LL^*(\overline{\nu}I_N-A^*)^{-1}C^*\\
&\hspace{-155mm}=(W-C({\lambda}I_N-A)^{-1}L)(W^*-L^*(\overline{\nu}I_N-A^*)^{-1}C^*)
    \end{split}
  \]
  and hence the result.
\end{proof}

Toward Corollary \ref{negH} we recall the following definition, which originates with the work of Krein; see \cite{ikl}:

\begin{definition}
  Let $\kappa\in\mathbb N_0$. A $\mathbb C^{p\times p}$-valued function $K({\lambda},{\nu})$ defined for ${\lambda},w$ in a set $\Omega$ is said to have $\kappa$ negative squares if it is Hermitian:
  \[
    K({\lambda},{\nu})=K(\nu,{\lambda})^*,\quad {\lambda},\nu\in\Omega,
    \]
    and if  for every $N\in\mathbb N$, $\xi_1,\ldots,\xi_N\in\mathbb C^p$ and $\nu_1,\ldots, \nu_N\in\Omega$, the $N\times N$ Hermitian matrix $(\xi^*_jK(\nu_j,\nu_k)
    \xi_k)_{j,k=1}^N$ has at most $\kappa$ strictly
    negative eigenvalues, and exactly $\kappa$ strictly
    negative eigenvalues for some choice of $N,\xi_1,\ldots, \xi_N,\nu_1,\ldots, \nu_N$.
  \end{definition}

Note that the case $\kappa=0$ corresponds to the classical notion of positive definite function.

  \begin{corollary}
    \label{negH}
  In the above notation, it follows from \eqref{decomp123} that the kernel function $\frac{\varphi({\lambda})+(\varphi({\nu}))^*}{1-{\lambda}\overline{\nu}}$ has at most $\kappa$ negative squares, where $\kappa$ is the number of strictly negative eigenvalues of
  $H$.
  \end{corollary}

In the terminology of Section \ref{krpon}, $\varphi$ is a generalized Carath\'eodory function when the number of negative squares of the kernel is strictly greater than $0$, and
  a Carath\'eodory function where it is equal to $0$.\smallskip

Finally we recall the following result, and give its proof for completeness. We refer also to \cite{gk1} for more details and information.

\begin{theorem}
  Let $R$ be a $\mathbb C^{p\times p}$-valued rational function, and assume that $R$ has no singularities on the unit circle.
  Let $R({\lambda})=D+C({\lambda}I_N-A)^{-1}B$ be a minimal realization of $R$ and let $P$ denote the Riesz projection corresponding to the spectrum of $A$ inside the open unit disk:
  \begin{equation}
    P=\frac{1}{2\pi i}\int_{|\lambda|=1}(\lambda I_N-A)^{-1}d\lambda.
  \end{equation}
  Then the restriction of $A$ to the range of $I_N-P$ is invertible, and $R({\lambda})=\sum_{u\in\mathbb Z}R_u\lambda^u$ with
  \begin{equation}
    R_u=\begin{cases}\, \delta_{u0}D-C((I_N-P)A|_{{\rm ran}\,I_N-P})^{-u-1}(I_N-P)B,\quad\, u\ge 0,\\
      \,C(PA)^{-u-1}PB,\quad u< 0.
    \end{cases}
    \label{coeff-fourier}
  \end{equation}
  (where $\delta_{uv}$ denotes Kronecker's symbol) and $\sum_{u\in\mathbb Z} \|R_u\|<\infty$.
  \label{Arche-de-la-defense}
\end{theorem}

\begin{proof}
  Let $P$ be the Riesz projection corresponding to the spectrum of $A$ inside the open unit disk. We write (with an abuse of notation
  for $(I_N-P)A|_{{\rm ran}\,I_N-P})^{-1}(I_N-P))$)
  \[
    \begin{split}
      (\lambda I_N-A)^{-1}&=(\lambda I_N-A)^{-1}(I_N-P+P)\\
      &=(\lambda I_N-PA)^{-1}P+(\lambda I_N-(I_N-P)A)^{-1}(I_N-P)\\
      &=\lambda^{-1}(I_N-\lambda^{-1}PA)^{-1}P-(I_N-\lambda ((I_N-P)A)^{-1})^{-1}((I_N-P)A)^{-1}(I_N-P)\\
      &=\sum_{u=0}^\infty \lambda^{-u-1}(PA)^uP-\sum_{u=0}^\infty \lambda^u((I_N-P)A)^{-u-1}(I_N-P).
      \end{split}
    \]
  \end{proof}

We now study minimal realizations of rational functions strictly positive on the unit circle without poles or zeros there.
As is known (see \cite{MR2663312,ot21}, and see \cite{MR0355675} for the general Wiener algebra setting), such functions admit Wiener-Hopf
factorizations, meaning in the rational case that one can write
\begin{equation}
  \label{spectral-facto-123}
  R(\lambda)=w(\lambda)(w(1/\overline{\lambda}))^*
\end{equation}
  where $w(\lambda)$ and its inverse are analytic in the exterior of the closed
  unit disk. See also \cite[Lemma 1.2, p. 145]{MR1294714} where the realization is centered at the origin. The function $w(\lambda)$ is called
  the left spectral factor, and is unique up to multiplication on the right by a unitary matrix.

\begin{proposition}
Let
\begin{equation}
  \label{bruxelles}
w(\lambda)=d+c(\lambda I_\ell-a)^{-1}b
\end{equation}
be  a minimal realization of  the  left spectral factor of $R(\lambda)$ in
\eqref{spectral-facto-123}. Then $\sigma(a)$ and $\sigma(a-bd^{-1}c)$ lie inside the open unit disk, and
\begin{equation}
  R(\lambda)=D+C(\lambda I_{2\ell}-A)^{-1}B,
\end{equation}
where
\begin{align}
  \label{prod111}
  A&=\begin{pmatrix}a&-bb^*a^{-*}\\0&a^{-*}\end{pmatrix}\\
  \label{prod222}
B&=\begin{pmatrix}b(d^*-b^*a^{-*}c^*)\\ a^{-*}c^*\end{pmatrix}\\
   C&=\begin{pmatrix}c&-db^*a^{-*}\end{pmatrix}\\
  D&=d(d^*-b^*a^{-*}c^*).
     \label{prod333}
\end{align}
\label{amsterdam}
\end{proposition}
\begin{proof}
Using the formula \eqref{formprod} for the realization
of a product, we have:
  \[
    \begin{split}
      R(\lambda)&=w(\lambda)(w(1/\overline{\lambda}))^*\\
&=\left(d+c(\lambda I_\ell-a)^{-1}b\right)\left(d+c\left(({1}/{\overline{\lambda}})I_\ell-a\right)^{-1}b\right)^*\\
      &=      \left(d+c(\lambda I_{\ell}-a)^{-1}b\right)\left(d^*+\lambda b^*\left(I_{\ell}-\lambda a^*\right)^{-1}c^*\right)\\
      &=      \left(d+c(\lambda I_{\ell}-a)^{-1}b\right)\left(d^*+\lambda b^*\left(a^{-*}-\lambda I_{\ell}\right)^{-1}a^{-*}c^*\right)\\
      &=      \left(d+c(\lambda I_{\ell}-a)^{-1}b\right)\left(d^*+b^*\left(\lambda I_{\ell}-a^{-*}+a^{-*}\right)\left(a^{-*}-\lambda I_{\ell}\right)^{-1}a^{-*}c^*\right)\\
      &=      \left(d+c(\lambda I_{\ell}-a)^{-1}b\right)\left(d^*-b^*a^{-*}c^*-b^*a^{-*}\left(\lambda I_{\ell}-a^{-*}\right)^{-1}a^{-*}c^*\right)\\
      &=D+C(\lambda I_{2\ell}-A)^{-1}B
\end{split}
\]
where $A,B,C$ and $D$ are given by \eqref{prod111}-\eqref{prod333}.
\end{proof}
\begin{proposition}
The Riesz projection $P$ is equal to
\begin{equation}
\begin{pmatrix}I_{\ell}&X\\0&0\end{pmatrix}
\end{equation}
where $X$ is the unique solution to the Stein equation
\begin{equation}
  \label{stein-equa}
  X-aXa^*=bb^*.
\end{equation}
\end{proposition}

\begin{proof}
  We have
  \[
  \begin{split}
    P&=\frac{1}{2\pi i}\int_{|\lambda|=1}\left(\lambda I_{2\ell}-A\right)^{-1}d\lambda\\
    &=\frac{1}{2\pi i}\int_{|\lambda|=1}\begin{pmatrix}(\lambda I_\ell-a)^{-1}&(\lambda I_{\ell}-a)^{-1}bb^*a^{-*}(\lambda I_\ell-a^{-*})^{-1}\\0&(\lambda I_\ell-a^{-*})^{-1}\end{pmatrix}d\lambda\\
    &=\begin{pmatrix}I_{\ell}&X\\0&0\end{pmatrix}
    \end{split}
  \]
  where
  \[
    \begin{split}
      X&=-\frac{1}{2\pi i}\int_{|\lambda|=1}(\lambda I_{\ell}-a)^{-1}bb^*a^{-*}(\lambda I_{\ell}-a^{-*})^{-1}d\lambda\\
      &=\frac{1}{2\pi i}\int_{|\lambda|=1}(\lambda I_{\ell}-a)^{-1}bb^*(I_{\ell}-\lambda a^{*})^{-1}d\lambda.
      \end{split}
\]
\end{proof}

We now give formulas for the Fourier coefficients in terms of $a,b,c$ and $d$, which are used in particular in Proposition \ref{prop5-3}.
It is useful to first rewrite $R(\lambda)$ in a more symmetric form

\begin{lemma}
In the notation of the Proposition \ref{amsterdam} we have:
  \begin{equation}
    R(\lambda)=\left(d+\sum_{u=1}^\infty \lambda^{-u}ca^{u-1}b\right)\left(d+\sum_{u=1}^\infty \lambda^{-u}ca^{u-1}b\right)^*.
    \end{equation}
  \end{lemma}

  \begin{proof}
We have for $|\lambda|=1$,
\[
  \begin{split}
    d^*-b^*a^{-*}c^*-b^*a^{-*}\left(\lambda I_{\ell}-a^{-*}\right)^{-1}a^{-*}c^*&=d^*-b^*a^{-*}c^*+b^*\left(I_{\ell}-\lambda a^{*}\right)^{-1}a^{-*}c^*\\
    &=d^*-b^*a^{-*}c^*+\sum_{u=0}^\infty \lambda^ub^*a^{*(u-1)}c^*\\
    &=d^*+\sum_{u=1}^\infty \lambda^ub^*a^{*(u-1)}c^*,
  \end{split}
  \]
  and hence the result since
    \[
      R(\lambda)=\left(d+c(\lambda I_{\ell}-a)^{-1}b\right)\left(d^*-b^*a^{-*}c^*-b^*a^{-*}\left(\lambda I_{\ell}-a^{-*}\right)^{-1}a^{-*}c^*\right).
      \]
\end{proof}
    Set
    \[
R(\lambda)=\sum_{k\in\mathbb Z}R_k\lambda^{k}
      \]
    where the $R_k$ are the Fourier coefficients of $R(\lambda)$. We have:

\begin{proposition}
  The Fourier coefficients of the spectral function  are given by
  \begin{align}
    \label{r0}
  R_0&=dd^*+cXc^*\\
  R_k&=(db^*+cXa^*)a^{*(k-1)}c^*,\quad k=1,2,\ldots\\
  R_{-k}&=R_k^*,\quad k=1,2,\ldots
            \label{r-k}
  \end{align}
  in terms of a realization \eqref{bruxelles} of the left spectral factor.
\end{proposition}

\begin{proof}
  By identifying the coefficient of $\lambda^k$ and the fact that $X=\sum_{t=0}^\infty a^tbb^*a^{*t}$. Indeed,
  \[
R_0=dd^*+\sum_{u=1}^\infty c a^{u-1}bb^*a^{*(u-1)}c=dd^*+cXc^*
\]
and, for $k>0$
\[
\begin{split}
  R_k&=db^*a^{*(k-1)}c^*+\sum_{u=1}^\infty ca^{(u-1)}bb^*a^{*(u+k-1)}c^*\\
  &=db^*a^{*(k-1)}c^*+cXa^{*k}c^*\\
    &=(db^*+cXa^*)a^{*(k-1)}c^*.
\end{split}
\]
By symmetry we have $R_k=R_k^*$ for $k\in\mathbb Z$.
\end{proof}

The case $d=0$ leads to simplifications, as illustrated in the following proposition.

\begin{proposition}
In the previous notations, we have
  \begin{equation}
    \label{dnot0}
    \varphi(\lambda)    =d(d^*-b^*a^{-*}c^*)+(db^*a^{-*}+cX)(I_\ell+\lambda a^*)
    (I_\ell-\lambda a^*)^{-1}c^*
\end{equation}
\begin{equation}
  \label{form02222022b}
  \varphi(\lambda)=-db^*a^{-*}c^*-cXc^*+2(db^*a^{-*}+cX)a^{-*}(a^{-*}-\lambda I_\ell)^{-1}c^*.
  \end{equation}
  Assume $d=0$. Then,
  \begin{equation}
    \label{form02222022a}
\varphi(\lambda)=cX(a^{-*}+\lambda I_\ell)(a^{-*}-\lambda I_\ell)^{-1}c^*
  \end{equation}
  rewritten as
  \begin{equation}
    \label{form02222022bb}
    \varphi(\lambda)=-cXc^*-2cXa^{-*}(\lambda I_\ell -a^{-*})^{-1}c^*
  \end{equation}
  corresponds to the choice $P=a^*X^{-1}a$ in the Kalman-Yakubovich-Popov conditions
  \eqref{kyp1}-\eqref{kyp3}.
\end{proposition}

\begin{proof} We can write $\varphi(\lambda)$ using the Fourier coefficients in the
  previous proposition as
  \[
  \begin{split}
    \varphi(\lambda)&=R_0+2\sum_{k=1}^\infty\lambda^kR_k\\
        &=dd^*+cXc^*+2\sum_{k=1}^\infty \lambda^k db^*a^{*(k-1)}c^*+2\sum_{k=1}^\infty\lambda^kcXa^{*k}c^*\\
        &=dd^*+2db^*a^{-*}(\lambda a^*)(I_\ell-\lambda a^*)^{-1}c^*+cXc^*+2 cX(\lambda a^*)(I_\ell-\lambda a^*)^{-1}c^*\\
        &=dd^*-db^*a^{-*}c^*+db^*a^{-*}\left(I_\ell+2(\lambda a^*)(I_\ell-\lambda a^*)^{-1}\right)c^*+cX(I_\ell+2(\lambda a^*)(I_\ell-\lambda a^*)^{-1}c^*\\
        &=d(d^*-b^*a^{-*}c^*)+(db^*a^{-*}+cX)(I_\ell+\lambda a^*)(I_\ell-\lambda a^*)^{-1}c^*
  \end{split}
  \]
  from where we have \eqref{form02222022a}. To prove \eqref{form02222022b} we write
  \[
  \begin{split}
    \varphi(\lambda)&=\varphi(\infty)+\varphi(\lambda)-\varphi(\infty)\\
    &=-2(db^*a^{-*}c^*+cXc^*)+2(db^*a^{-*}+cX)(a^{-*}+\lambda I_\ell)(a^{-*}-\lambda I_\ell )^{-1}c^*+2(db^*a^{-*}+cX)c^*\\
   &= -db^*a^{-*}c^*-cXc^*+2(db^*a^{-*}+cX)a^{-*}(a^{-*}-\lambda I_\ell)^{-1}c^*.
    \end{split}
  \]
  We now assume $d=0$ and check that conditions \eqref{kyp1}-\eqref{kyp3} hold with the above realization
\[
  (a^{-*},-\sqrt{2}c^*,\sqrt{2}cXa^{-*},-cXc^*)
\]
of $\varphi$ and $P=aX^{-1}a^*$, with $W=0$.
We begin with \eqref{kyp3}. With $W=0$ we need to check that
\[
0=-2cXc^*+2cXa^{-*}a^*X^{-1}a^{-1}aXc^*,
  \]
  which holds. We now check that \eqref{kyp2} holds

  \[
  -\sqrt{2}a^{-*}a^*X^{-1}aa^{-1}Xc=-\sqrt{2}c^*
  \]
  which holds. We now check that $P-a^{-*}Pa^{-1}\le 0$, which will give the existence of an $L$
  in \eqref{kyp1}. Since $X^{-1}$ exists,
  \[
X\ge aXa^*
\]
is equivalent to $I_N\ge X^{-1/2}aX^{1/2}X^{1/2}a^*X^{-1/2}$, and so,
since the adjoint of a contraction is a contraction
\[
X^{1/2}a^{*}X^{-1}aX^{1/2}\le I_N
\]
or
\[
a^*X^{-1}a\le X^{-1}.
\]
But this implies
\[
P-a^{-*}Pa^{-1}=a^*X^{-1}a-X^{-1}\le 0.
\]

\end{proof}
\begin{remark}
{\rm The above computations for the realizations of $\varphi$ are redone in a different way later in the paper; see Section \ref{gohgoh}.}
\end{remark}

\begin{remark}{\rm
Taking a unitary dilation of $a^{-*}$ lead to an expression of the kind considered by Dijksma, Langer and de Snoo.
  }
  \end{remark}
\subsection{Rational discrete analytic functions}
\label{sub-real}
As mentioned in the introduction the following result, set here in the setting of matrices and convergent power series, plays a central role (see \cite[(3.20), p. 352]{MR450576}):

\begin{theorem}
  Let $f(m,n)$, $m,n=0,1,\ldots$ be a family of $p\times q$ matrices with complex entries, indexed by $\mathbb N_0\times\mathbb N_0$, and such that the power series
  $k_f(\lambda,\mu)=\sum_{m,n=0}^\infty f(m,n)\lambda^m\overline{\mu}^n$ converges in a neighborhood of the origin of $\mathbb C^2$. Then $f(m,n)$ defines a $\mathbb C^{p\times q}$-valued
  discrete analytic function in $\Lambda_{++}$ if and only if there exist functions $\Phi_L(\lambda)$ and $\Phi_R(\lambda)$ analytic in a neighborhood of the origin, and such that
\begin{equation}
k_f(\lambda, \mu)=\frac{\Phi_L(\lambda)+\overline{\Phi_R(\mu)}}{1+i\lambda-i\overline{\mu}-\lambda\overline{\mu}}.
\end{equation}
The functions $\Phi_L$ and $\Phi_R$ are given by \eqref{q1} and \eqref{q2}, and convergence of the power series \eqref{q1} and \eqref{q2} implies the convergence of the generating function
in a neighborhood of the origin.
\end{theorem}

\begin{proof} The arguments in \cite{MR450576} show that there is convergence of the generating function when $\Phi_L$ and $\Phi_R$ converge in a neighborhood of the origin.\end{proof}

We now recall a number of background definitions, and send the reader to \cite{MR0078441,MR0013411} for more information. A path (say in $\Lambda_{++}$) is a sequence of points $\gamma=(z_0,\ldots, z_n)$ in
$\Lambda_{++}$ with $|z_k-z_{k-1}|=1$ for $k=1,\ldots, n$. The discrete integral along the path $\gamma$ is defined by
\begin{equation}
\int_\gamma fdz=\sum_{k=1}^n\frac{f(z_{k-1})+f(z_k)}{2}(z_k-z_{k-1}).
  \end{equation}
  A new operator introduced in \cite{alpay2021discrete} is
  \begin{equation}
(Zf)(z)=\frac{f(0)-f(z)}{2}+\int_0^zf\delta s
\end{equation}

Pointwise product of two discrete analytic functions need not be discrete analytic, and we will use the product \eqref{proddaf} below:

\begin{definition} (see \cite{alpay2021discrete})
  Let $u\in\mathbb N$. We define:
  \begin{align}
    Z^u1&=z^{(u)}\\
    (Z^{(u)}1)\odot f)&=Z^uf.
                        \label{proddaf}
  \end{align}
\end{definition}
The following results and definition are taken from \cite{alpay2021discrete,alpay2021discrete2}.
\begin{proposition}
We have, with $z=m+ni$,
\begin{equation}
  \label{totoche321}
\sum_{u=0}^\infty   t^uz^{(u)}=(1+t)^m\left(\frac{1+\alpha_+t}{1+\alpha_-t}\right)^n,
\end{equation}
where
\[
\alpha_+=\frac{1+i}{2}\quad and\quad \alpha_-=\frac{1-i}{2}.
\]
\end{proposition}

\begin{proposition}
  \label{mnmnmn}
  Let $z=m+in$ and let $A\in\mathbb C^{N\times N}$ be such that $\left\{-2\alpha_\pm\right\}\cap \sigma(A)=\emptyset$.
  Then, it holds that
\begin{equation}
  (I_N-zA)^{-\odot}=(I_N+A)^m(I_N+\alpha_+A)^n(I_N+\alpha_-A)^{-n}
\end{equation}
and the matrix-valued function $(I_N-zA)^{-\odot}$ is discrete analytic in
\(
\Lambda_+=\left\{z\in\Lambda\,;\, {\rm Re}\, z\ge 0\right\}
\)  
(as opposed to $\Lambda_{++}$).
\end{proposition}

Motivated by the classical realization theory of rational functions, discrete analytic rational functions were introduced in \cite{alpay2021discrete2} in terms of realization similar to \eqref{real-0} as follows:

\begin{definition} A matrix-valued discrete analytic function will be called rational if it can be written in the form
  \begin{equation}
    \label{real-00}
    f(z)=D+C(I_N-zA)^{-\odot}\odot (zB)
  \end{equation}
  where $A,B,C$ and $D$ are matrices of compatible sizes, and where $\left\{-2\alpha_\pm\right\}\cap \sigma(A)=\emptyset$.
  \end{definition}

In \cite{alpay2021discrete2} it has been proved that $\odot$-products and $\odot$-inverses of discrete analytic rational functions are still discrete analytic rational.

  \begin{proposition}
    Let $(C,A,B)\in\mathbb C^{p\times N}\times\mathbb C^{N\times N}\times \mathbb C^{N\times q}$.
The functions $z\mapsto (I_N-zA)^{-\odot}$ and $z\mapsto C(I_N-zA)^{-\odot}B$ are rational and discrete analytic.
  \end{proposition}

  \begin{proof}
    For the first function, we take $D=0$, $C=I_N$ and $B=A$ in \eqref{real-00}. We have
    \[
(I_N-zA)^{-\odot}\odot (zA)=(I_N-zA)^{-\odot}\odot (zA-I_N+I_N))=-I_N+(I_N-zA)^{-\odot}
\]
and hence $(I_N-zA)^{-\odot}$ is rational. The second claim follows by multiplying on the left and right by constant matrices.
    \end{proof}

    As a corollary of the above and of Proposition \ref{mnmnmn}:

    \begin{corollary}
Let $A\in\mathbb C^{N\times N}$ be such that $\left\{-2\alpha_\pm\right\}\cap \sigma(A)=\emptyset$.
 The $\mathbb C^{N\times N}$-valued function
    \begin{equation}
f(m,n)=(I_N+A)^m(I_N+\alpha_+A)^n(I_N+\alpha_-A)^{-n}
\label{ratio-daf}
\end{equation}
is rational discrete analytic, and so is the function
\[
  Cf(m,n)B=C(I_N+A)^m(I_N+\alpha_+A)^n(I_N+\alpha_-A)^{-n}B
\]
where $C\in\mathbb C^{c\times N}$ and $B\in\mathbb C^{N\times b}$.
\end{corollary}

To make the connection with \cite{alpay2021discrete,alpay2021discrete2} and get a more symmetric formula we set
\begin{align}
  A_1&=I_N+A,\\
  A_2&=(I+\alpha_{+}A)(I+\alpha_{-}A)^{-1}.
  \end{align}
Then, $A_1$ and $A_2$ satisfy \eqref{newconditioncom11} and \eqref{ratio-daf} becomes
\begin{equation}
  \label{ratio-daf-2}
  f(m,n)=A_1^mA_2^n,\quad m,n=0,1,\ldots
\end{equation}

\begin{lemma}
  Let $A_1=\sqrt{2}\alpha^*-iI_N$ and $A_2=\sqrt{2}\alpha+iI_N$ where $\alpha\in\mathbb C^{n\times n}$. Then, $A_1$ and $A_2$ satisfy \eqref{newconditioncom11} if and only if $\alpha$ is unitary.
\end{lemma}

\begin{proof} The assertion follows immediately from
  \begin{eqnarray}
\nonumber
    I_N+iA_1-iA_2-A_1A_2&=&I_N+i\sqrt{2}\alpha^*+I_N-i\sqrt{2}\alpha+I_N-(2\alpha\alpha^*+i\sqrt{2}\alpha^*-i\sqrt{2}\alpha+I_N)\\
                        &=&2(I_N-\alpha^*\alpha).
                              \label{equa567}
\end{eqnarray}
\end{proof}

\begin{remark}
  {\rm The right hand side in \eqref{equa567} becomes $0$ when considering the unitary dilation $U$ of $\alpha$ in a Krein space (see Theorem \ref{unit-dil} for the definition). Then,
  \begin{align}
      P(\sqrt{2}U^{-1}-iI_{\mathfrak H})^{m}|_{\mathbb C^N}&=(\sqrt{2}\alpha^*-iI_N)^m\\
      P(\sqrt{2}U+iI_{\mathfrak H})^{n}|_{\mathbb C^N}&=(\sqrt{2}\alpha+iI_N)^n.
  \end{align}
Then the pair $(A_1,A_2)$ is replaced by $(\sqrt{2}U^{-1}-iI_{\mathfrak H}, \sqrt{2} U+iI_{\mathfrak H})$.
  }
  \end{remark}


  For every $u\in\mathbb N$ the function $z^{(u)}$ is rational, but the function \eqref{real-00} will not have, in general, an expansion along the functions $z^{(u)}$, $u=1,2,\ldots$.
  In the next proposition, which explores this question, $r(A)$ denotes the spectral radius of the matrix $A$.

  \begin{proposition}
    In the above notation, assume $r(A)<\sqrt{2}$. Then,
    \begin{equation}
    \label{real-000}
    D+C(I_N-zA)^{-\odot}\odot (zB)=D+\sum_{u=1}^\infty z^{(u)}CA^{u-1}B.
    \end{equation}
 \end{proposition}

 \begin{proof}
   From \eqref{totoche321} one has that, for $z=m+in$ with $m>0$ and $n\not=0$
   \begin{equation}
\limsup_{u\rightarrow\infty}|z^{(u)}|^{1/u}=\frac{1}{\sqrt{2}},
\label{groningue7890}
\end{equation}
   and so the series $\sum_{u=0}^\infty z^{(u)}A^u$ converges for $r(A)<\sqrt{2}$ and diverges
   for $r(A)>\sqrt{2}$.
   \end{proof}

   We note that the bound \eqref{groningue7890} is uniform in $z=m+in$ with $m>0$ and $n\not=0$.\\


A discrete analytic function on $\Lambda_{++}$ is uniquely determined by its values $f(m,0)$ and $f(0,n)$, with $m,n\in\mathbb N_0$.
So:

\begin{theorem}
  \label{ex123}
  The function
  \[
f(m,n)=CA_1^mA_2^nB,\quad (m,n)\in\mathbb N_0\times \mathbb N_0,
\]
where $(C,A_1)$ is an observable pair of matrices and $(A_2,B)$ is a controllable pair of matrices defines a discrete analytic function if and only if \eqref{newconditioncom11} holds:
\begin{equation*}
I+iA_1-iA_2-A_1A_2=0.
\end{equation*}
Then $iI+A_1$ is invertible, $A_1$ and $A_2$ commute and
$$A_2=(iI+A_1)^{-1}(I+iA_1).$$
\end{theorem}

\begin{proof}
Since $f$ is a discrete analytic function, \eqref{def-daf} holds and so
    \[
CA_1^n\left(I+iA_1-iA_2-A_1A_2\right)A_2^mB=0,\quad m,n=0,1,\ldots
\]
The observability of the pair $(C,A_1)$ implies that
\[
  \left(I+iA_1-iA_2-A_1A_2\right)A_2^mB=0,\quad m=0,1,2,\ldots
\]
The controllability of the pair $(A_2,B)$ implies in turn that $I+iA_1-iA_2-A_1A_2=0$.\smallskip

We now show that $A_1$ and $A_2$ commute.
Assume $uA_1=-iu$, thus

$$u+u-iuA_2+iuA_2=0,$$
and hence $u=0$. Hence, $iI+A_1$ is invertible and we rewrite \eqref{newconditioncom11} as

$$I+iA_1=(iI+A_1)A_2$$
so $A_2=(iI+A_1)^{-1}(I+iA_1)$, and in particular $A_1$ and $A_2 $ commute.
\end{proof}

  \begin{corollary}
Assume that $A_1$ and $A_2$ satisfy \eqref{newconditioncom11}. They commute and hence
\[
(I+iA_1)(I-iA_2)=2A_1A_2
\]
or, provided inverses exist
\[
\sigma(A_1)\sigma(A_2^*)^*=I,
\]
where $\sigma(\lambda)$ is defined by \eqref{sigma}.
\end{corollary}

\begin{remark}{\rm In the special case
    where  $A_1=A_2^*$, i.e. $A_1$ is normal, we are back to equation \eqref{s-a-+c-t-u}.}
  \end{remark}

\begin{remark}{\rm
    Conversely, given the boundary values $f(m,0)=CA_1^mB$ and $f(0,n)=CA_2^nB$ with $n,m\in\mathbb N_0$, it
 is not clear if one can obtain a closed formula for $f(m,n)$ when  $A_1$ and $A_2$  do not satisfy \eqref{newconditioncom11}.}
  \end{remark}

\subsection{Krein spaces, Pontryagin spaces and generalized Carath\'eodory functions}
\label{krpon}
\begin{definition}
A $\mathbb C^{p\times p}$-valued function analytic in a neighborhood $\mathcal N$ of the origin is called a Carath\'eodory function if the kernel
    \begin{equation}
    \label{kphi}
\frac{\varphi({\lambda})+\varphi(\nu)^*}{1-{\lambda}\overline{\nu}}
\end{equation}
is positive definite in $\mathcal N$, and a generalized Carath\'eodory function if \eqref{kphi} has a finite number, say $\kappa$, of negative squares in $\mathcal N$.
\label{cara789}
\end{definition}

We first consider the Hilbert space case, i.e. $\kappa=0$. The analyticity hypothesis can be weakened in the definition and one has:

  \begin{theorem}
  Let $\varphi$ be a $\mathbb C^{p\times p}$-valued function {\sl defined} on a uniqueness set $\mathcal Z_\varphi$ of the open unit disk, and such that the kernel \eqref{kphi}
  is positive definite on $\mathcal Z_\varphi$. Then, $\varphi$ has a (unique) analytic extension to the open unit disk, for which the kernel \eqref{kphi} is still positive definite. This extension
  can be written as
\begin{equation}
  \label{herglotzform}
 \varphi({\lambda})=iX+\frac{1}{2\pi}\int_0^{2\pi}\frac{e^{it}+{\lambda}}{e^{it}-{\lambda}}dM(t),
 \end{equation}
 where $M$ is an increasing matrix-valued function on $[0,2\pi]$ of finite variation, the integral is a Stieltjes integral and $X$ is a Hermitian
 matrix.
 \label{belleville}
\end{theorem}

For Schur functions, the analytic extension follows from \cite{donoghue}. The claim for Carath\'eodory functions follows by Cayley transform
formula \eqref{herglotzform} is due to Herglotz (see \cite{herglotz}, \cite{hspnw}).\\

To go to the indefinite metric setting we first recall some definitions.

\begin{definition}
  \label{def-krein}
  A vector space $\mathfrak V$ endowed with a Hermitian form $[\cdot,\cdot]$ is a Krein space if it can be written as
  \begin{equation}
  \mathfrak V=  \mathfrak V_+\oplus  \mathfrak V_-
  \end{equation}
where the sum is direct and orthogonal, and where $(\mathfrak V_+,[\cdot,\cdot])$ and
$(\mathfrak V_-,-[\cdot,\cdot])$ are Hilbert spaces. It is called a Pontryagin space when the latter is moreover finite dimensional.
\end{definition}
We refer to  \cite{azih,bognar,ikl} (note that in \cite{ikl} the space $\mathfrak V_+$ is assumed to be finite dimensional).\\

Theorem \ref{belleville} has a counterpart in the Pontryagin setting, but an analyticity hypothesis
is needed. The analytic extension result will not hold when the kernel has a finite number of negative squares as the following example shows:
    \[
      \varphi({\lambda})=\begin{cases}\, 1,\,\,\,\,\, {\lambda}\in B(0,1)\setminus\left\{0\right\},\\
                              \, 0,\,\,\,\,\, {\lambda}=0.\end{cases}
      \]
      See \cite[p. 82]{adrs}, \cite[p. 701]{atv1}, \cite{MR2481906}.
      But one has (see \cite[Satz 2.1, p.  361]{kl1}):

      \begin{theorem}
        \label{kl-extension}
Let $\varphi$ be a generalized Carath\'eodory function. Then, it has a unique meromorphic extension to the open unit disk, and this extension can be written as
\begin{equation}
  \label{realcara}
    \varphi({\lambda})=iX+\frac{1}{2}C(U-{\lambda}I)^{-1}(U+{\lambda}I)C^*
  \end{equation}
  where $X=X^*\in\mathbb C^{p\times p}$, and $U$ is a unitary operator in a Pontryagin space $\mathfrak P$ of index $\kappa$, and where $C$ is a bounded operator from $\mathfrak P$ into $\mathbb C^p$.
  \end{theorem}

Note that in the previous theorem one can assume only continuity at the origin to ensure a meromorphic extension.
Note also that Herglotz integral formula admits generalizations to the case of  generalized Carath\'eodory functions, given by Krein and Langer; see for instance \cite[(4.11) p. 215]{MR47:7504},
but we will work with the realization formula \eqref{realcara}.\\

One can go one step further and consider the Krein space setting. Among the important differences between Pontryagin and Krein spaces, we mention that a unitary operator in a Pontryagin space
is continuous, while there exists unbounded unitary operators in Krein spaces; see
\cite{MR944858} a discussion of the latter. The following is the
specialization to matrix-valued functions of a result of Dijksma, Langer and de Snoo
\cite{MR89a:47055}, \cite[Theorem 1, p. 126]{MR903068}.

\begin{theorem}
  Let $\varphi$ be a $\mathbb C^{n\times n}$-valued function analytic in $|{\lambda}|<r_0$ with
  $r_0\in(0,1)$. Then there is a Krein space $\mathfrak K$, a bounded unitary operator $U$ in
  $\mathfrak K$ and a bounded operator $C$ from $\mathfrak K$ to $\mathbb C^n$ such that
  \begin{equation}
    \label{krein-987}
    \varphi({\lambda})=i{\rm Im}\, \varphi(0)+C(U+{\lambda}I)(U-{\lambda}I)^{-1}C^*.
  \end{equation}
  \label{dlds}
\end{theorem}

We will also need the notion of unitary dilation, which we now recall.

\begin{theorem} \cite[p. 78]{MR264438}
  Let $T$ be a closed densely defined operator in the Hilbert space $\mathfrak H$. There exists
  a Hilbert space $\mathfrak G$ containing $\mathfrak H$ and a closed densely defined operator $U$ in $\mathfrak G$ and a signature operator $J$ in $\mathfrak G$ such that:\\
  $(1)$ $\mathfrak G$ with the inner product induced by $J$ is a Krein space and
  $\mathfrak H$ is $J$-positive.\\
  $(2)$ $U$ is unitary and $U^{-1}$ is densely defined.\\
  $(3)$ $U$ is a dilation of $T$, meaning that
  \begin{align}
    T^n&=P_{\mathfrak H}U^n|_{\mathfrak H},\quad n=0,1,\ldots\\
    T^{*n}&=P_{\mathfrak H}U^{-n}|_{\mathfrak H},\quad n=,1,2,\ldots
    \end{align}
  $(4)$ It holds that
  \[
\bigvee\left\{ {\rm ran}\, U^n\mathfrak H,\,\, n\in\mathbb Z\right\}=\mathfrak G.
\]
\label{unit-dil}
  \end{theorem}

The operator $U$ is called a unitary dilation of $T$.
When $T$ is a contraction the space $\mathfrak G$ is a Hilbert space. See \cite[p. 16]{nf}.
In the Hilbert space setting $U$ is uniquely determined (up to a Hilbert space isomorphism) by the conditions in the theorem, but not in the general case; see the discussion \cite[p. 85]{MR264438}.
Of importance here is the following consequence of the construction given in \cite[p. 78]{MR264438}.

\begin{corollary}
Assume that $T$ in the above theorem is everywhere defined and bounded. Then the unitary dilation is also bounded.
 \end{corollary}

  \begin{proof}
    We repeat the construction of \cite{MR264438}, specialized to a bounded operator $T$, and check that the extension is indeed bounded. With $E(t)$ and $F(t)$ the spectral resolutions associated
    to the positive operators $(TT^*)^{1/2}$ and $(T^*T)^{1/2}$ respectively, one first introduces the operators (see \cite[\S 2]{MR264438})
    \begin{align}
      J_T&=\int_0^\infty {\rm sgn}\,(1-t^2)dE(t)\\
      Q_T&=\int_0^\infty \sqrt{|1-t^2|}dE(t)\\
      Q_{T^*}&=\int_0^\infty {\rm sgn}\,(1-t^2)dF(t).
    \end{align}
    Now these operators are bounded since $(TT^*)^{1/2}$ and $(T^*T)^{1/2}$ are bounded, and hence have spectral resolutions with compact support in $[0,\infty)$. The operator $U$ in \cite{MR264438} is an
    extension of the operator
    \[
      V=\begin{pmatrix}
        \ddots&  & & & & & & &\\
          \ddots    & 0& & & & & & & \\
              &I&0&  & & & & & \\
              & &I&0     & & & & & \\
              & & &Q_{T^*}&\boxed{T}   & & & & \\
              & & &-J_TT^*&Q_T&0& & & \\
              & & &       &  &I&0 & &\\
              & & &       &  & &I &0 &\\
              & & &       &  & & &\ddots &\ddots
        \end{pmatrix}
      \]
      which is bounded since $T$ is bounded.
    \end{proof}

    \subsection{Structured matrices, displacement rank and dualities}
    We discuss connections with the topics in the title of this subsection and with reproducing kernel spaces of pairs.
  Assume $f(m,n)\in\mathbb C^{p\times p}$ and
  \[
F_N=(f(m,n))_{m,n=0}^N.
    \]

  Set $Z\in\mathbb C^{(N+1)p\times (N+1)p}$ to be the backward shift matrix
  \begin{equation}
    Z=\begin{pmatrix}0&I_p&0&\cdots&0\\
      0&0&I_p&\cdots&0\\
      \vdots&\vdots&\vdots & &\vdots\\
      0&0&0&\cdots&I_p\\
      0&0&0&\cdots &0\end{pmatrix}.
    \end{equation}

    In the Hermitian case we have (see \cite[p. 148]{MR1197502}):
\begin{proposition}
    In the Hermitian case we have
  \begin{equation}
    \label{equa-daf}
a(Z^*)F_Na(Z^*)^*-b(Z^*)F_Nb(Z^*)^*=V^*JV
\end{equation}
or, equivalently, equation \eqref{equa-daf} can be rewritten as
         \begin{equation}
           \label{equa-daf2}
F_N+iZ^*F_N-iF_NZ-Z^*F_NZ=V^*JV
\end{equation}
    with
    \[
      J=\begin{pmatrix}0&I_p\\I_p&0\end{pmatrix}\quad{\rm and}\quad
      V=\begin{pmatrix} \frac{f(0,0)}{2}&f(0,1)-if(0,0)&\cdots &f(0,N)-if(0,N-1)\\
          -I_p&   0_{p\times p}&\cdots&0_{p\times p}\end{pmatrix}.
         \]
       \end{proposition}

       \begin{proof} As is well known in matrix theory, the action of $Z^*$ on the left and $Z$ on the right are given by
         \[
           Z^*F_N=\begin{pmatrix}0&0&\cdots &0\\
             f(0,0)&f(0,1)&\cdots&f(0,N)\\
            \vdots &\vdots & &\vdots\\
             f(N-1,0)&f(N-1,1)&\cdots&f(N-1,N)\end{pmatrix}
         \]

and
                  \[
                    F_NZ=\begin{pmatrix}0&f(0,0)&f(0,1)&\cdots&f(0,N-1)\\
            \vdots & & &\vdots \\
             0&f(N,0)&f(N,1)&\cdots&f(N,N-1)\end{pmatrix}.
         \]
         Thus
         \[
           Z^*F_NZ=\begin{pmatrix}0&0&\cdots &0\\
             0&f(0,0)&\cdots &f(0,N-1)\\
             \vdots&\vdots &\cdots &\vdots\\
             0&f(N-1,0)&\cdots &f(N-1,N-1)\end{pmatrix}.
         \]
         Assuming $f(m,n)$ discrete analytic we get
         \[
           \begin{split}
             F_N+iZ^*F_N-iF_NZ-Z^*F_NZ&=\\
             &\hspace{-4cm}=-\begin{pmatrix}f(0,0)&f(0,1)-if(0,0)&\cdots &f(0,N)-if(0,N-1)\\
  f(1,0)+if(0,0)&0&\cdots &0\\
  \vdots&\vdots &\cdots &\vdots\\
  f(N,0)+if(N-1,0)&0&\cdots &0
\end{pmatrix}\\
&=V^*JV
\end{split}
\]
with $V$ and $J$ as in the theorem.
         \end{proof}

         \begin{theorem} (see \cite[Theorem 5.3, p.  150]{ad-laa4})
           Assume that $F_N$ is Hermitian and invertible. Then,
           \begin{equation}
V(I_N+i\lambda I_N-iZ-\lambda Z)^{-1}F_N^{-1}(I_N+i\nu I_N-iZ-\nu Z)^{-*}V^*=\frac{J-\Theta(\lambda)J\Theta(\nu)^*}{1+i\lambda-i\overline{\nu}-\lambda\overline{\nu}}
\end{equation}
with
\begin{equation}
  \Theta(\lambda)=I_{2p}-(1+iz-i\overline{\mu_0}-z\overline{\mu_0})V(I_N+i\lambda I_N-iZ-\lambda Z)^{-1}F_N^{-1}(I_N+i\mu_0 I_N-iZ-\mu_0 Z)^{-*}V^*
\end{equation}
where $\mu_0$ is any point in $C(-i,\sqrt{2})$ such that $I_N+i\mu_0 I_N-iZ-\mu_0 Z$ is invertible. The choice $\mu_0=1$ gives
\begin{equation}
  \Theta(\lambda)=I_{2p}-(1-\lambda)V(I_N+i\lambda I_N-iZ-\lambda Z)^{-1}F_N^{-1}(I_N-Z)^{-*}V^*.
  \end{equation}
           \end{theorem}

           \begin{remark} {\rm Note that from $V$ one can recover the whole matrix $F_N$ using the discrete analytic Cauchy equations.}
             \end{remark}

Reproducing kernel pairs of spaces in duality (and in particular reproducing kernel spaces) and structured matrices are two closely related topics. Important examples of the connections occur with
the Gohberg-Semencul  \cite{MR0353038} and Gohberg-Heinig (noncommutative algebra version of the latter; see
\cite{GH}) inversions formula versus the Christoffel-Darboux formula for orthogonal polynomials. Similar connections exist for Hankel matrices and more generally for structured matrices.
\\

\begin{lemma}
  \label{lemmalrlr}
  Assume that \eqref{richelieu-drouot} converges in a neighborhood of the closure of
  $B(-i,\sqrt{2})$. Then,
  \begin{equation}
\label{symm-789}
    \Phi_L({\lambda})+(\Phi_R({\lambda}^\prime))^*=0,
  \end{equation}
  where ${\lambda}^\prime =\frac{1-i\overline{{\lambda}}}{\overline{{\lambda}}-i}$.
\end{lemma}

\begin{proof}
  This follows from
  \[
    \sigma({\lambda})\overline{\sigma({\lambda}^\prime)}=1
  \]
  and \eqref{opera-bastille111111}.
  \end{proof}

We note
that the elements invariant under the map ${\lambda}\mapsto {\lambda}^\prime$ are exactly the points of $C(-i,\sqrt{2})$, and that $i^\prime=0$. Furthermore,
${\lambda}^{\prime\prime}={\lambda}$ for ${\lambda}\not=-i$ and, on the Riemann sphere, $-i^\prime=\infty$ and $\infty^\prime=-i$.\\

The following results are the specialization of \cite[Theorem 4.1, p. 37]{ad9} to $a({\lambda})$ and $b({\lambda})$ given by \eqref{opera-garnier}; see also \cite[Theorem 5.2, p. 24]{ad-jfa} and the discussion page 26 in that paper, and \cite[Theorem 4.1 p. 424]{ad-laa4}.\smallskip

In the statement of Theorem \ref{thmmm}, the matrices
  $(V_\ell,A_\ell,B_\ell)\in\mathbb C^{p\times q}\times\mathbb C^{q\times q}\times\mathbb C^{q\times q}$, $\ell=L,R$
  satisfy,
\begin{equation}
\det\,((1+i{\lambda})A_\ell-\sqrt{2}{\lambda}B_\ell)\not\equiv 0,\quad \ell=L,R,
 \end{equation}
and, with $c\in\mathbb C^q$,
  \begin{equation}
    V_\ell((1+i{\lambda})A_\ell-\sqrt{2}{\lambda}B_\ell)^{-1}c\equiv 0\,\,\Longrightarrow\,\, c=0,\quad \ell=L,R.
  \end{equation}
  Furthermore we set
  \begin{equation}
    F_\ell({\lambda})=V_\ell((1+i{\lambda})A_\ell-\sqrt{2}{\lambda}B_\ell)^{-1},\quad \ell=L,R,
  \end{equation}
  and denote by $\mathfrak M_\ell$ the linear span of the columns of $F_\ell$ for $\ell=L,R$.
  \begin{theorem}
    \label{thmmm}
        Let $P\in\mathbb C^{q\times q}$ be an invertible matrix, and endow the space $\mathfrak M_L\times\mathfrak M_R$ with the sesquilinear form
    \begin{equation}
      [F_Lc,F_Rd]=d^*Pc,\quad c,d\in\mathbb C^q.
    \end{equation}
    Then $\mathfrak M_L\times\mathfrak M_R$ is a reproducing kernel pair space with reproducing kernel pair
    \begin{equation}
      K_L({\lambda},{\nu})=F_L({\lambda})P^{-1}F_R({\nu})^*\quad and\quad K_R({\lambda},{\nu})=F_R({\lambda})P^{-1}F_L({\nu})^*,
    \end{equation}
    meaning that
    \begin{align}
      [F_L(\cdot)\xi,K_R(\cdot,{\nu})\eta]&=\eta^*F_L({\nu})\xi \\
     [ K_L(\cdot,{\nu})\eta,F_R(\cdot,{\nu})\xi ]&=(F_R({\nu})\xi)^*\eta.
      \end{align}
\end{theorem}

  \begin{theorem}
    Let $J\in\mathbb C^{p\times p}$ be a signature matrix.
In the notation of the preceding theorem, the reproducing kernel has the form
    \begin{eqnarray}
      \label{rk-ell}
      K_L({\lambda},{\nu})&=&\frac{J-W_L({\lambda})JW_R({\nu})^*}{1+i{\lambda}-i\overline{\nu}-{\lambda}\overline{\nu}},\\
      K_R({\lambda},{\nu})&=&\frac{J-W_R({\lambda})JW_L({\nu})^*}{1+i{\lambda}-i\overline{\nu}-{\lambda}\overline{\nu}},
    \end{eqnarray}
    for some $\mathbb C^{p\times p}$-valued functions $W_L({\lambda})$ and $W_R({\lambda})$
    if and only if $P$ is a solution of the matrix equation
    \begin{equation}
      \label{tuileries}
      A_R^*PA_L-B_R^*PB_R=V_R^*JV_L.
    \end{equation}
    The functions $W_L$ and $W_R$ are then uniquely given by the formulas
\begin{eqnarray}
  W_L({\lambda})&=&(I_p-(1+i{\lambda}-i\overline{\mu}-{\lambda}\overline{\mu})F_L({\lambda})P^{-1}F_R(\mu)^*J)C_L,\\
            W_R({\lambda})&=&(I_p-(1+i{\lambda}-i\overline{\mu}-{\lambda}\overline{\mu})F_L({\lambda})P^{-1}F_L(\mu)^*J)C_R
                      \end{eqnarray}
                      where $\mu$ belongs to the circle $|{\lambda}+i|=\sqrt{2}$ and is a point of analyticity of $F_L$ and $F_R$, and
where $U_L$ and $U_R$ belong to $\mathbb C^{p\times p}$ and satisfy $U_LJU_R^*=J$.
\end{theorem}

The choice $\mu=1$ gives $1+i{\lambda}-i\overline{\mu}-{\lambda}\overline{\mu}=(1-i)(1-{\lambda})$ and then,

\begin{eqnarray}
  W_L({\lambda})&=&(I_p-(1-i)(1-{\lambda})F_L({\lambda})P^{-1}F_R(i)^*J)C_L,\\
            W_R({\lambda})&=&(I_p-(1-i)(1-{\lambda})F_L({\lambda})P^{-1}F_L(i)^*J)C_R.
                      \end{eqnarray}
Important cases are when
\begin{equation}
  A_\ell=1-i
  \mathsf A_\ell\quad{\rm and}\quad B_\ell=\sqrt{2}\mathsf A_\ell,
  \end{equation}
  for some $\mathsf A_\ell\in\mathbb C^{q\times q}$, $\ell=L,R$.
  Then,
  \begin{equation}
(1+i{\lambda})(1-i\mathsf A_\ell)-\sqrt{2}{\lambda}\mathsf A_ \ell=1+i{\lambda}-i\mathsf A_\ell-{\lambda}\mathsf A_\ell,\quad \ell=L,R,
  \end{equation}
  and the matrix equation \eqref{tuileries} becomes
  \begin{equation}
P+i\mathsf A_R^*P-iP\mathsf A_\ell-\mathsf A_R^*P\mathsf A_L=V_R^*JV_L.
\end{equation}

  \begin{lemma}
    Let ${\lambda}\in\mathbb C\setminus\left\{-i\right\}$, not a pole of $W_L$ or $W_R$. Then
    \begin{equation}
W_L({\lambda})J(W_R({\lambda}^\prime))^*=J.
    \end{equation}
    Assume that $A_L$ and $B_R$ invertible. Then, $W_L$ has no pole and is invertible at $0$, while
    $W_R$ has no pole and is invertible at $i$.\\
    Assume that $A_R$ and $B_L$ invertible. Then, $W_L$ has no pole and is invertible at $i$,
    while
    $W_R$ has no pole and is invertible at $0$.\\
    \end{lemma}

    \begin{proof}
      The first claim follows from \eqref{rk-ell}, rewritten as
      \[
W_L({\lambda})JW_R({\nu})^*=J-({1+i{\lambda}-i\overline{\nu}-{\lambda}\overline{\nu}}) K_L({\lambda},{\nu}),
      \]
      with $w={\lambda}^\prime$.\smallskip

      Using the first claim the second claim follows with ${\lambda}=0$ since we then have:
      \begin{equation}
        W_L(0)JW_R(i)^*=J.
        \end{equation}
      \end{proof}


We now take
\begin{equation}
  J=\begin{pmatrix}0&1\\1&0\end{pmatrix}\quad{\rm and}\quad W_\ell=\begin{pmatrix}a_\ell&b_\ell\\ c_\ell&d_\ell\end{pmatrix},\,\,\ell=L,R.
\end{equation}
Mutliplying \eqref{rk-ell} on the left by a $J$-neutral vector $e=\begin{pmatrix}u&v\end{pmatrix}
\in\mathbb C^2$ (i.e. $e^*Je=0$)
we get

\[
  \begin{split}
    eK_L({\lambda},{\nu})e^*&=\frac{-eW_L({\lambda})JW_R({\nu})^*e^*}{1+i{\lambda}-i\overline{\nu}-{\lambda}\overline{\nu}}\\
    &=\frac{-
      (\overline{u}a_L({\lambda})+\overline{v}c_L({\lambda}))
      ({u}\overline{b_R({\nu})}+v\overline{d_R({\nu})})-
            (\overline{u}b_L({\lambda})+\overline{v}d_L({\lambda}))
      ({u}\overline{a_R({\nu})}+v\overline{c_R({\nu})})}
    {1+i{\lambda}-i\overline{\nu}-{\lambda}\overline{\nu}}\\
    &=
      (\overline{u}a_L({\lambda})+\overline{v}c_L({\lambda}))
    \frac{\Phi_L({\lambda})+\overline{\Phi_R({\nu})}}{1+i{\lambda}-i\overline{\nu}-{\lambda}\overline{\nu}}
        ({u}\overline{a_R({\nu})}+v\overline{c_R({\nu})})
    \end{split}
  \]
  with
  \begin{equation}
    \Phi_m({\lambda})=-\frac{
\overline{u}b_\ell({\lambda})+\overline{v}d_\ell({\lambda})}
    {\overline{u}a_\ell({\lambda})+\overline{v}c_\ell({\lambda})},\quad \ell=L,R.
  \end{equation}
  provided
  \[
\overline{u}a_\ell({\lambda})+\overline{v}c_\ell({\lambda}))\not\equiv0,\quad \ell=L,R.
  \]
The case where these functions do not vanish at $0$ leads to a family of examples, as explained in the next proposition:

  \begin{proposition}
    Assume the four matrices $A_\ell,B_\ell$, $\ell=L,R$, to be invertible. Then there exists
    a $J$-neutral vector $e$ such that the function
  \[
  \frac{e^*K_L({\lambda},{\nu})e}{(\overline{u}a_L({\lambda})+\overline{v}c_L({\lambda}))(
{u}\overline{a_R({\nu})}+{v}\overline{c_R({\nu})})}= \frac{\Phi_L({\lambda})+\overline{\Phi_R({\nu})}}{1+i{\lambda}-i\overline{\nu}-{\lambda}\overline{\nu}}
     \]
     can be developed in powers series in the variable ${\lambda}$ and $\overline{\nu}$ near $(0,0)$, and hence is the generating function of a discrete analytic function in $\Lambda_{++}$.
  \end{proposition}

 \begin{proof}
   By Lemma \ref{lemmalrlr} we have
   \begin{eqnarray}
     W_L(0)JW_R(i)^*&=&J
     \label{LR1}\\
     W_L(i)JW_R(0)^*&=&J.
     \label{LR2}
     \end{eqnarray}
         By \eqref{LR1} we cannot have both $a_L(0)=0$ and $c_L(0)=0$. Similarly, by
         \eqref{LR2}, we cannot have both $a_R(0)=c_R(0)=0$. In fact we can have at most
         two of the four numbers $a_L(0),a_R(0),c_L(0),c_R(0)$ equal to $0$.
         To find a vector $e$ we consider the following cases:\\

         {Case 1:} $a_L(0)=a_R(0)=0$.  We take $e=\begin{pmatrix}0\\1         \end{pmatrix}$.\smallskip

         {Case 2:} $c_L(0)=c_R(0)=0$.  We take $e=\begin{pmatrix}1\\0
         \end{pmatrix}$.\smallskip

         {Case 3:} $a_L(0)=c_R(0)=0$ or $a_R(0)=c_L(0)=0$.
         We take $e=\begin{pmatrix}1\\1\end{pmatrix}$.\smallskip

         {Case 4:} $a_L(0)=0$ but $a_R(0)c_R(0)\not=0$, or
         $a_R(0)=0$ but $a_L(0)c_L(0)\not=0$.
         We take $e=\begin{pmatrix}0\\1         \end{pmatrix}$.\smallskip

         {Case 5:} $c_L(0)=0$ but $a_R(0)c_R(0)\not=0$, or
         $c_R(0)=0$ but $a_R(0)c_R(0)\not=0$.
         We take $e=\begin{pmatrix}1\\0         \end{pmatrix}$.\smallskip

         \end{proof}

\section{Realizations for coefficients of discrete analytic functions}
\setcounter{equation}{0}

We first consider the non symmetric case and then the Krein and Pontryagin spaces respectively.

\subsection{Spaces in dualities}
For spaces in duality we refer in particular to the papers \cite{aron1,zbMATH03043485}. For the
following definition, see \cite[Definition 1, p. 1245]{aro3}, where it is given in the setting
of Banach spaces.

\begin{definition}
  \label{groningen}
  Let $\mathfrak H_L$ and $\mathfrak H_R$ be two Hilbert spaces of $\mathbb C^p$-valued functions defined on a set $\Omega$, with norms $\|\cdot\|_L$ and $\|\cdot\|_R$, and
  let $[\cdot,\cdot]$ be a sesquilinear form on $\mathfrak H_L\times \mathfrak H_R$. Then, the pair $(\mathfrak H_L\times \mathfrak H_R,[\cdot,\cdot])$ is a reproducing kernel Hilbert space of pairs
  if there exists a pair $(K^L,K^R)$ of $\mathbb C^{n\times n}$-valued functions  defined on $\Omega\times\Omega$ and such that:
    \begin{align}
      [f,K_R(\cdot,{\nu})\eta]&=\eta^*f({\nu}) \\
     [ K_L(\cdot,{\nu})\eta,g ]&=g({\nu})^*\eta
      \end{align}
for every $f\in\mathfrak H_L$, every $g\in\mathfrak H_R$, every $w\in\Omega$ and every $\eta\in\mathbb C^p$.
\end{definition}

We consider now the general case where $k_f({\lambda},{\nu})$ is given by \eqref{opera-bastille111}. The same computations hold for the proofs of \eqref{formulemoments}.
With $X=0$ in \eqref{q1} we have near the origin
\[
  \sum_{m=0}^\infty f(m,0){\lambda}^m
  =\frac{\Phi_L({\lambda})+\frac{f(0,0)}{2}}{1+i{\lambda}}
\]

\begin{theorem}
  \label{tmtmtm}
  There exist Krein spaces $\mathfrak K_j$, bounded unitary operators $U_j\,:\,\mathfrak K_j\,\,\rightarrow\,\, \mathfrak K_j$ and bounded operators $C_j\,:\, \mathfrak K_j\,\,\mapsto \,\,
  \mathbb C^p$ ($j=L,R$) such that
  \begin{align}
\label{q1q1q1}
    f(m,0)&=C_L(\sqrt{2}U_L^{-1}-iI)^mC_L^*,\quad m=0,1,\ldots\\
    f(0,n)&=C_R(\sqrt{2}U_R+iI)^nC_R^*,\quad\,\,\,\, n=0,1,\ldots
  \end{align}
  \end{theorem}

  \begin{proof}
Setting ${\lambda}=0$ we obtain
\[
f(0,0)=\Phi_L(0)+\frac{f(0,0)}{2}.
\]
Then, we can write in a neighborhood of the origin
\[
  \begin{split}
    \sum_{m=0}^\infty f(m,0){\lambda}^m&=
  \frac{\Phi_L({\lambda})+\frac{f(0,0)}{2}}{1+i{\lambda}}\\&=\frac{f(0,0)+\sum_{m=1}^\infty \Phi_{L,m} {\lambda}^m}{1+i{\lambda}}\\
  &=\left(f(0,0)+\sum_{m=1}^\infty \Phi_{L,m} {\lambda}^m\right)\left(\sum_{m=0}^\infty (-i)^m{\lambda}^m\right)\\
  &=\sum_{m=0}^\infty {\lambda}^m\left(f(0,0)(-i)^m+\sum_{p=1}^m \Phi_{L,p}(-i)^{m-p}\right).
  \end{split}
\]
Let now
\[
\Phi_L(\sigma^{-1}({\lambda}))=iM+C_L(U_L+{\lambda}I)(U_L-{\lambda}I)^{-1}C_L^*
\]
be a realization of $\varphi_L$ obtained from Dijksma-Langer-de Snoo theorem (see Theorem \ref{kl-extension} above).
We use this representation with $M=0$. The computations are the same as for the proofs of \eqref{formulemoments}. We have:
\begin{equation}
  \label{nantes}
  \begin{split}
    \Phi_L({\lambda})&=\frac{1}{2}C_L\left(U_L-\frac{\sqrt{2}{{\lambda}}}{1+i{\lambda}}I\right)^{-1}\left(U_L+\frac{\sqrt{2}{{\lambda}}}{1+i{\lambda}}I\right)C_L^*\\
        &=\frac{1}{2}C_L\left((1+i{\lambda})U_L-\sqrt{2}{{\lambda}}I\right)^{-1}\left((1+i{\lambda})U_L+\sqrt{2}{{\lambda}}I\right)C_L^*\\
        &=\frac{1}{2}C_L(U_L+{\lambda}(iU_L-\sqrt{2})I)^{-1}(U_L+{\lambda}(iU_L+\sqrt{2})I)C_L^*\\
        &=\frac{1}{2}C_LC_L^*+\frac{1}{2}C_L\left((U_L+{\lambda}(iU_L-\sqrt{2}I))^{-1}(U_L+{\lambda}(iU_L+\sqrt{2})I)-I\right)C_L^*\\
        &=\frac{1}{2}C_LC_L^*+\sqrt{2}{\lambda}C_L(U_L+{\lambda}(iU_L-\sqrt{2}I))^{-1}C_L^*\\
        &=\frac{1}{2}C_LC_L^*+\sqrt{2}{\lambda}C_LU_L^{-1}(I-{\lambda}(\sqrt{2}U_L^{-1}-iI))^{-1}C_L^*\\
        &=\frac{1}{2}C_LC_L^*+\sum_{n=0}^\infty {\lambda}^{n+1}\sqrt{2}C_LU_L^{-1}(\sqrt{2}U_L^{-1}-iI)^nC_L^*\\
        &=\sum_{n=0}^\infty\Phi_{L,n}{\lambda}^n.
    \end{split}
  \end{equation}

  By \eqref{nantes}

  \[
    \begin{split}
      f(m,0)&=(-i)^mf(0,0)+\sum_{p=1}^m \Phi_{L,p}(-i)^{m-p}\\
      &=      f(0,0)(-i)^m+\sum_{p=1}^m(-i)^{m-p}\sqrt{2}C_LU_L^{-1}(\sqrt{2}U_L^{-1}-iI)^{p-1}C_L^*\\
      &=      f(0,0)(-i)^m+(-i)^m\sum_{p=1}^m\sqrt{2}C_LU_L^{-1}(\sqrt{2}U_L^{-1}-iI)^{-1}(\sqrt{2}U_L^{-1}-iI)^{p}(-i)^pC_L^*\\
      &=      f(0,0)(-i)^m+(-i)^m\sum_{p=1}^m\sqrt{2}C_LU_L^{-1}(\sqrt{2}U_L^{-1}-iI)^{-1}(i\sqrt{2}U_L^{-1}+I)^{p}C_L^*\\
      &=      f(0,0)(-i)^m+\\
      &\hspace{-1cm}+(-i)^m\sqrt{2}C_LU_L^{-1}(\sqrt{2}U_L^{-1}-iI)^{-1}(i\sqrt{2}U_L^{-1}+I)(I-(i\sqrt{2}U_L^{-1}+I)^m)(I-(i\sqrt{2}U_L^{-1}+I)^{-1}C_L^*\\
      &=      f(0,0)(-i)^m-(-i)^mC_L(I-(i\sqrt{2}U_L^{-1}+I)^mC_L^*,
      \end{split}
    \]
    and so
    \begin{equation}
      f(m,0)=C_L(\sqrt{2}U_L^{-1}-iI)^mC_L^*,\quad m=0,1,\ldots
      \label{newform}
      \end{equation}
\end{proof}

\subsection{The Krein space setting}

We now consider the case of functions $f$ such that the kernel $k_f$ is Hermitian.
\begin{lemma}
  Assume that the kernel $k_f({\lambda},\nu)$ is Hermitian in a uniqueness set $\mathcal Z_f\subset B(-i,\sqrt{2})$, meaning that
  \begin{equation}
    \label{sym67}
    k_f({\lambda},\nu)=k_f(\nu,{\lambda})^*,\quad {\lambda},w\in \mathcal Z_f.
  \end{equation}
  Then, one can choose $\Phi_L$ and $\Phi_R$ such that $\Phi_L({\lambda})=\Phi_R({\lambda})$ for ${\lambda}\in\mathcal Z_f$. Moreover, $f(0,0)=f(0,0)^*$ and $X=0$ in \eqref{q1}-\eqref{q2}.
  \end{lemma}

\begin{proof}
  From \eqref{sym67} and \eqref{opera-bastille} we get
  \begin{equation}
\Phi_L({\lambda})+\Phi_R(\nu)^*=\Phi_R({\lambda})+\Phi_L(\nu)^*,\quad {\lambda},w\in\mathcal Z_f,
\end{equation}
and so
\begin{equation}
  \Phi_L({\lambda})-\Phi_R({\lambda})=(\Phi_L(\nu)-\Phi_R(\nu))^*,\quad {\lambda},\nu\in\mathcal Z_f.
\end{equation}
It follows that $\Phi_L({\lambda})-\Phi_R({\lambda})=M$, where $M=M^*$ is a constant. Replacing $\Phi_L({\lambda})$ by $\Phi_L({\lambda})-M/2$ and $\Phi_R({\lambda})$ by $\Phi_R({\lambda})+M/2$ we obtain the result. Finally
writing $\Phi_L(0)=\Phi_R(0)$ in \eqref{q1}-\eqref{q2} we obtain
\[
\frac{f(0,0)}{2}+iX=\frac{f(0,0)^*}{2}+iX^*,
\]
and so
\[
\frac{f(0,0)-f(0,0)^*}{2}=i(X^*-X),
\]
which implies that $X=0$ and $f(0,0)=f(0,0)^*$.
  \end{proof}

In the Hermitian case the kernel $k_f({\lambda},\nu)$ can thus be rewritten as
\begin{equation}
  \label{opera-bastille-2}
  k_f({\lambda},\nu)=\frac{\Phi({\lambda})+\Phi(\nu)^*}{1+i{\lambda}-i\overline{\nu}-
    {\lambda}\overline{\nu}},
\end{equation}
with
\begin{equation}
\label{phipos}
\Phi({\lambda})=(1+i{\lambda})\left(\sum_{m=0}^\infty f(m,0){\lambda}^m\right)-\frac{f(0,0)}{2}
\end{equation}
with $f(0,0)=f(0,0)^*$.\\

In a way analogous to \eqref{nantes} we now write $\Phi$ as in \eqref{phipos} with $q=0$ to get near the origin
\[
  \sum_{m=0}^\infty f(m,0){\lambda}^m
  =\frac{\Phi({\lambda})+\frac{f(0,0)}{2}}{1+i{\lambda}}.
\]
Setting ${\lambda}=0$ we obtain
\[
f(0,0)=M_0+\frac{f(0,0)}{2}.
\]

Then, still near $0$, and since
\[
\Phi(0)+\frac{f(0,0)}{2}=f(0,0),
\]
we can write
\[
  \begin{split}
    \sum_{m=0}^\infty f(m,0){\lambda}^m&=
  \frac{\Phi({\lambda})+\frac{f(0,0)}{2}}{1+i{\lambda}}\\&=\frac{f(0,0)+\sum_{m=1}^\infty \Phi_m {\lambda}^m}{1+i{\lambda}}\\
  &=\left(f(0,0)+\sum_{m=1}^\infty \Phi_m {\lambda}^m\right)\left(\sum_{m=0}^\infty (-i)^m{\lambda}^m\right)\\
  &=\sum_{m=0}^\infty {\lambda}^m\left(f(0,0)(-i)^m+\sum_{p=1}^m \Phi_{p}(-i)^{m-p}\right).
  \end{split}
\]
Hence
\begin{eqnarray}
\label{nantes-2}
  f(m,0)&=&f(0,0)(-i)^m+\sum_{p=1}^m \Phi_{p}(-i)^{m-p},\quad m=1,2,\ldots,
\\
\nonumber
\end{eqnarray}
(compare with Proposition \ref{5-12}, Proposition \ref{propfin} and Remark \ref{remfin}).

  \begin{definition}
    \label{symdef}
The $\mathbb C^{p\times p}$-valued function discrete analytic in $\Lambda_{++}$ will be called symmetric if
\[
  f(m,0)=f(0,m)^*,\quad m=0,1,\ldots
  \]
\end{definition}
Symmetry means that in formula \eqref{opera-bastille} it holds that $\Phi_L({\lambda})=\Phi_R({\lambda})=\Phi(\lambda)$ in a neighborhood of the origin.
We apply the realization result of Dijksma, Langer and de Snoo (see Theorem \ref{dlds}) and obtain:
  \begin{theorem}
    Let $f$ be a symmetric discrete analytic function in $\Lambda_{++}$ such that the generating function $k_f({\lambda},\nu)$ converges in a neighborhood of $(0,0)$. Then there is a bounded unitary operator in a Krein
    space $\mathfrak K$ and a bounded operator from $\mathfrak K$ into $\mathbb C^p$ such that \eqref{gare-du-nord1} holds:
    \begin{equation*}
      \label{gare-du-nord}
      f(m,n)=C(\sqrt{2}U^{-1}-iI)^m(\sqrt{2}U+iI)^nC^*,\quad n,m=0,1,\ldots
    \end{equation*}
  \end{theorem}

  \begin{proof}
    We start from a realization \eqref{gare-de-l-est} of $\Phi(\sigma^{-1}({\lambda}))$:
    \[
      \Phi(\sigma^{-1}({\lambda}))=iX+\frac{1}{2}C(U+{\lambda}I)(U-{\lambda}I)^{-1}C^*.
      \]
      The same computations as in Theorem \ref{tmtmtm} gives now
  \begin{align}
\label{q1q1q1q1}
    f(m,0)&=C\sqrt{2}U^{-1}-iI)^mC^*,\quad m=0,1,\ldots\\
    f(0,n)&=C(\sqrt{2}U+iI)^nC^*,\quad\,\,\,\, n=0,1,\ldots
    \end{align}
    \end{proof}

Here too, the existence of a discrete analytic extension depends on the spectrum of $U$. We note that
\begin{equation}
  \label{gener-fn-1}
  \sum_{m=0}^\infty {\lambda}^mf(m,0)=C(I-{\lambda}(\sqrt{2}U^*-iI))^{-1}C^*.
\end{equation}

\begin{theorem}
  \label{montparnasse12}
  The generating function can be written as
  \begin{equation}
    \label{handaye1}
k_f({\lambda},\nu)=C(I-{\lambda}(\sqrt{2}U^*-iI))^{-1}(I-{\overline{\nu}}(\sqrt{2}U+iI))^{-1}C^*.
    \end{equation}
\end{theorem}

\begin{proof}
  From \eqref{gener-fn-1} we have
  \[
    \begin{split}
      k_f({\lambda},\nu)&=\frac{C\left((1+i{\lambda})(I-{\lambda}(\sqrt{2}U-iI))^{-1}+(1-i\overline{\nu})(I-\overline{\nu}(\sqrt{2}U^*+iI))^{-1}-1\right)C^*}{1+i{\lambda}-i\overline{\nu}-{\lambda}\overline{\nu}}\\
      &=\frac{C(1-{\lambda}(\sqrt{2}U-iI))^{-1}N({\lambda},w,U)(1-\overline{\nu}(\sqrt{2}U^*+i))^{-1}C^*}{(1+i{\lambda}-i\overline{\nu}-{\lambda}\overline{\nu})}
      \end{split}
    \]
    with
    \[
\begin{split}
  N({\lambda},w,U)&=(1+i{\lambda})(I-\overline{\nu}(\sqrt{2}U^*+iI))+(I-i\overline{\nu})(I-{\lambda}(\sqrt{2}U-iI))-\\
  &\hspace{5mm}-(I-{\lambda}(\sqrt{2}U-iI))(I-\overline{\nu}(\sqrt{2}U^*+iI))\\
  &=I+i{\lambda}I-\overline{\nu}\sqrt{2}U^*-i\ow I-i{\lambda}\overline{\nu}\sqrt{2}U^*-{\lambda}\ow I+\\
  &\hspace{5mm}+I-i\overline{\nu}I-{\lambda}\sqrt{2}U+i{\lambda}\ow\sqrt{2}U+{\lambda}\ow+i{\lambda}I-\\
 &\hspace{5mm}-I+\ow\sqrt{2}U^*+i\ow I+{\lambda}\sqrt{2}U-i{\lambda} I-\\
  &\hspace{5mm}-{\lambda}\ow(2I+i\sqrt{2}U-i\sqrt{2}U^*+I)\\
&=  (1+i{\lambda}-i\overline{\nu}-{\lambda}\overline{\nu})I.
  \end{split}
      \]
  \end{proof}

To make connection with \cite{alpay2021discrete,alpay2021discrete2}, we remark the following.
Let $U$ be unitary in a Krein space and let $A=\sqrt{2}U^*-iI-I$. Then,
  \[
    \begin{split}
      (I+\alpha_+A)(I+\alpha_-A)^{-1}
      &=\left(I+\frac{1+i}{2}(\sqrt{2}U^*-(i+1)I\right)\left(I+\frac{1-i}{2}\sqrt{2}U^*-(i+1)I\right)^{-1}\\
      &=\left((2-(1+i)^2)I+(1+i)\sqrt{2}U^*\right)\left((2-(1+i)(1-i))I+(1-i)\sqrt{2}U^*\right)^{-1}\\
      &=\left(2(1-i)I+(1+i)\sqrt{2}U^*\right)\left((1-i)\sqrt{2}U^*\right)^{-1}\\
&=\sqrt{2}U+iI.
      \end{split}
    \]
    So
    \begin{equation}
      (\sqrt{2}U^*-iI)^{m}(\sqrt{2}U+i)^n=(I+A)^m\left((I+\alpha_+A)(I+\alpha_-A)^{-1}\right)^n.
      \end{equation}

      \subsection{Pontryagin space case}
      We now specialize the Krein space setting to the Pontryagin setting.
\begin{theorem}
  Let $f$ be a discrete $\mathbb C^{p\times p}$-valued analytic function for which the kernel \eqref{opera-bastille} has a finite number of negative squares in a neighborhood of the origin in
  $B(-i,\sqrt{2})$. Then there is a generalized Carath\'eodory function $\varphi$ with the same number of negative squares, and such that
  \begin{equation}
k_f({\lambda},{\nu})=\frac{\varphi(\sigma({\lambda}))+\varphi(\sigma({\nu}))^*}{1+i{\lambda}-i\overline{\nu}-{\lambda}\overline{\nu}},
\end{equation}
where ${\lambda},w$ varies in $B(-i,\sqrt{2})$, from which are removed the points with image a pole of $\varphi$.
 \end{theorem}

 \begin{proof}
   With $\Phi$ as in
   Let $g({\lambda})=\Phi(\sigma^{-1}({\lambda}))$. We rewrite \eqref{opera-bastille-2} as
   \[
a({\lambda})k_f({\lambda},{\nu})\overline{a({\nu})}=\frac{g(\sigma({\lambda}))+g(\sigma({\nu}))^*}{1-\sigma({\lambda})\overline{\sigma({\nu})}}.
\]
The function $g$ is analytic in a neighborhood, say $U\subset B(0,1)$, of the origin and the kernel $\frac{g({\lambda})+g({\nu})^*}{1-{\lambda}\overline{\nu}}$ has a finite number of negative squares in $U$. It follows from
Theorem \ref{kl-extension} that $g$ has a meromorphic extension to $B(0,1)$, with the same number of negative squares. Denoting by $\varphi$ this extension we get the result.
\end{proof}

Using the realization formula \eqref{realcara} we have the following. The proof is similar to previous computations and omitted.

 \begin{theorem}
   \label{th4-3}
   We have
   \begin{equation}
\Phi({\lambda})=\sum_{m=0}^\infty \Phi_{m}{\lambda}^{m},
\end{equation}
with $\Phi_0=\frac{1}{2}CC^*$ and
\begin{equation}
  \label{phi-n-p}
  \Phi_m=CU^*(\sqrt{2}U^*-iI)^{m-1}C^*,\quad m=1,2,\ldots
\end{equation}
 \end{theorem}

   \section{The positive case}
   \setcounter{equation}{0}

\subsection{Formula for the coefficients}

In the positive case, that is when the generating function $k_f(\lambda,\nu)$ converges and is positive definite in a neighborhood of the origin, using Herglotz integral formula \eqref{herglotzform} we can prove the next result:

 \begin{theorem}
Let $\Phi$ be given in the form \eqref{herglotzform}. It holds that
   \begin{equation}
\Phi({\lambda})=\sum_{m=0}^\infty \Phi_{m}{\lambda}^{m},
\end{equation}
with $\Phi_0=\frac{1}{2\pi}\int_0^{2\pi}dM(t)$ and
\begin{equation}
  \label{phi-n}
  \Phi_m=\frac{\sqrt{2}}{\pi}\int_0^{2\pi}e^{-it}(\sqrt{2}e^{-it}-i)^{m-1}dM(t),\quad m=1,2,\ldots
\end{equation}
 \end{theorem}

 \begin{proof} We write
 \[
 \begin{split}
\Phi({\lambda})&=iX+\frac{1}{2\pi}\int_0^{2\pi}\frac{(1+i{\lambda})e^{it}+\sqrt{2}{\lambda}}{(1+i{\lambda})e^{it}-\sqrt{2}{\lambda}}dM(t)\\
&=iX+\frac{1}{2\pi}\int_0^{2\pi}\left(\frac{e^{it}(1+i{\lambda})-\sqrt{2}{\lambda}+2\sqrt{2}{\lambda}}{(1+i{\lambda})e^{it}-\sqrt{2}{\lambda}}\right)dM(t)\\
&=iX+\frac{1}{2\pi}\int_0^{2\pi}\left(1+\frac{2\sqrt{2}{\lambda}}{e^{it}}\frac{1}{1-{\lambda}(\sqrt{2}e^{-it}-i)}\right)dM(t)\\
&=iX+M_0+\frac{\sqrt{2}}{\pi}\sum_{m=0}^\infty {\lambda}^{m+1}\int_0^{2\pi}e^{-it}(\sqrt{2}e^{-it}-i)^mdM(t)\\
&=iX+M_0+\sum_{m=0}^\infty \Phi_{m+1}{\lambda}^{m+1},
 \end{split}
\]
with $X=X^*\in\mathbb C^{p\times p}$,
the matrices $\Phi_m$ given by \eqref{phi-n} and $M_0=\frac{1}{2\pi}\int_0^{2\pi}dM(t)$.
\end{proof}

\begin{remark}
  {\rm We note that \eqref{phi-n} is a specialization of the general formula \eqref{phi-n-p} to the positive case.}
  \end{remark}

The  $\Phi_m$, and hence the $f(m,0)$,  are functions of the trigonometric moments $M_0,\ldots, M_m$ for $m=1,2\ldots$.
\begin{theorem}
  Assuming the representation \eqref{herglotzform} for $\Phi$, the boundary value of the discrete analytic functions are given by
\begin{equation}
    \label{marseille}
f(m,0)=\frac{1}{\pi}\int_0^{2\pi}(\sqrt{2}e^{-it}-i)^mdM(t),\quad m=0,1,\ldots
\end{equation}
\end{theorem}

\begin{proof}
We first compute
\begin{equation}
  \label{f2m0}
  f(0,0)=2M_0=\frac{1}{\pi}\int_0^{2\pi}dM(t).
\end{equation}

Then we have:

\[
\begin{split}
  \sum_{\ell=1}^m \Phi_{\ell}(-i)^{m-\ell}&=\frac{(-i)^m}{\pi}e^{-it}\left(\sum_{\ell=1}^m(\sqrt{2}e^{-it}-i)^{\ell-1}(-i)^{-\ell}\right)dM(t)\\
  &=\frac{(-i)^m}{\pi}\int_0^{2\pi}\sqrt{2}e^{-it}(\sqrt{2}e^{-it}-i)^{-1}\left(\sum_{\ell=1}^m\frac{(\sqrt{2}e^{-it}-i)^\ell}{(-i)^\ell}\right)dM(t)\\
  &=\frac{(-i)^m}{\pi}\int_0^{2\pi}\sqrt{2}e^{-it}(\sqrt{2}e^{-it}-i)^{-1}\left(\sum_{\ell=1}^m(1+i\sqrt{2}e^{-it})^\ell\right)dM(t)\\
  &=\frac{(-i)^m}{\pi}\int_0^{2\pi}\sqrt{2}e^{-it}(\sqrt{2}e^{-it}-i)^{-1}(1+i\sqrt{2}e^{-it})\frac{1-(1+i\sqrt{2}e^{-it})^m}{1-(1+i\sqrt{2}e^{-it})}dM(t)\\
  &=-\frac{(-i)^m}{\pi}\int_0^{2\pi}\left(1-(1+i\sqrt{2}e^{-it})^m\right)dM(t).
\end{split}
\]

So, using \eqref{nantes-2} and \eqref{f2m0}, we have

\begin{equation}
  \begin{split}
    f(m,0)&=(-i)^m\left(f(0,0) -\frac{1}{\pi}\int_0^{2\pi}\left(1-(1+i\sqrt{2}e^{-it})^m\right)dM(t)\right)\\
    &=\frac{(-i)^m}{\pi}\int_0^{2\pi}(1+i\sqrt{2}e^{-it})^mdM(t),
  \end{split}
\end{equation}
and so we get \eqref{marseille}.
\end{proof}

We thus get to a generalized moment problem of the kind studied by Krein and Nudelman; see \cite[Theorem 1.1 p. 58 and Theorem 1.2 p. 59]{knud}.\\

We note that
\begin{equation}
  \label{gener-fn}
\sum_{m=0}^\infty \lambda^mf(m,0)=\frac{1}{\pi}\int_0^{2\pi}\frac{dM(t)}{1-\lambda(\sqrt{2}e^{-it}-i)},
\end{equation}
which can also be rewritten as
\begin{equation}
  \frac{1}{\pi} \int_0^{2\pi}\frac{dM(t)}{1+i\lambda-\sqrt{2}\lambda e^{-it}}\quad{\rm or}\quad  \frac{1}{\pi(1+i\lambda)}
  \int_0^{2\pi}\frac{dM(t)}{1-\sigma(\lambda) e^{-it}}.
  \end{equation}
\begin{theorem}
  \label{montparnasse}
It holds that
  \begin{equation}
    \label{handaye}
    k_f({\lambda},{\nu})=\frac{1}{\pi}    \int_0^{2\pi}\frac{dM(t)}{(1-{\lambda}(\sqrt{2}e^{-it}-i))(1-\overline{\nu}(\sqrt{2}e^{it}+i))}
    \end{equation}
and the $f(m,n)$ are given by \eqref{formulemoments}, namely:
  \begin{equation*}
f(m,n)=\frac{1}{\pi}\int_0^{2\pi}(\sqrt{2}e^{-it}-i)^m(\sqrt{2}e^{it}+i)^ndM(t),\quad m,n=0,1,\ldots
  \end{equation*}
\end{theorem}

\begin{proof}
  From \eqref{gener-fn}, \eqref{opera-bastille-2} and \eqref{phipos} we have
  \[
    \begin{split}
      k_f({\lambda},{\nu})&=\dfrac{1}{\pi}\frac{
        \displaystyle{\int_0^{2\pi}}\left(\dfrac{1+i{\lambda}}{1-{\lambda}(\sqrt{2}e^{-it}-i)}+\dfrac{1-i\overline{\nu}}{1-\overline{\nu}(\sqrt{2}e^{it}+i)}-1\right)dM(t)}{1+i{\lambda}-i\overline{\nu}-{\lambda}\overline{\nu}}\\
      &\\
      &=\frac{1}{\pi}
\dfrac{\displaystyle{\int_0^{2\pi}}\dfrac{N({\lambda},w,t)dM(t)}{(1-{\lambda}(\sqrt{2}e^{-it}-i))(1-\overline{\nu}(\sqrt{2}e^{it}+i))}}{(1+i{\lambda}-i\overline{\nu}-{\lambda}\overline{\nu})}
      \end{split}
    \]
    with
    \[
\begin{split}
  N({\lambda},w,t)&=(1+i{\lambda})(1-\overline{\nu}(\sqrt{2}e^{it}+i))+(1-i\overline{\nu})(1-{\lambda}(\sqrt{2}e^{-it}-i))-\\
  &\hspace{5mm}-(1-{\lambda}(\sqrt{2}e^{-it}-i))(1-\overline{\nu}(\sqrt{2}e^{it}+i))\\
  &=1+i{\lambda}-\overline{\nu}\sqrt{2}e^{it}-i\ow-i{\lambda}\overline{\nu}\sqrt{2}e^{it}-{\lambda}\ow+\\
  &\hspace{5mm}+1-i\overline{\nu}-{\lambda}\sqrt{2}e^{-it}+i{\lambda}\ow\sqrt{2}e^{-it}+{\lambda}\ow+i{\lambda}-\\
  &\hspace{5mm}-1+\ow\sqrt{2}e^{it}+i\ow+{\lambda}\sqrt{2}e^{-it}-i{\lambda}-\\
  &\hspace{5mm}-{\lambda}\ow(2+i\sqrt{2}e^{-it}-i\sqrt{2}e^{it}+1)\\
&=  1+i{\lambda}-i\overline{\nu}-{\lambda}\overline{\nu}.
  \end{split}
      \]
  \end{proof}

  \begin{remark}{\rm
 Although a consequence of the previous theorem we now check directly that the family
  $(f(m,n))_{m,n\in\mathbb N_0}$ satisfies \eqref{def-daf}. With
  $K_1=\sqrt{2}e^{-it}-i$ and $K_2=\sqrt{2}e^{it}+i$ (note that $K_2=\overline{K_1}$),
  this amounts to check that
    \[
1+iK_1-iK_2-K_1K_2=0,
\]
 which clearly holds since:
\[
  \begin{split}
    1+iK_1-iK_2-K_1K_2&=1+i(\sqrt{2}e^{-it}-i)-i(\sqrt{2}e^{it}+i)-(\sqrt{2}e^{-it}-i)(\sqrt{2}e^{it}+i)\\
    &=1+i\sqrt{2}e^{-it}+1-i\sqrt{2}e^{it}+1-2-i\sqrt{2}e^{-it}+i\sqrt{2}e^{it}-1\\
    &=0.
  \end{split}
\]
}
\end{remark}

\begin{corollary}
  In the notation and hypothesis of Theorem \ref{montparnasse} it holds that
  \begin{equation}
2 M_0 \left(\sqrt{2}-1\right)^{2n}\le f(n,n)\le 2 M_0(\sqrt{2}+1)^{2n},\quad n=0,1,\ldots
    \end{equation}
Moreover
  \begin{equation}
|f(m,n)|\le 2\pi M_0 (\sqrt{2}+1)^{m+n},\quad n=0,1,\ldots
    \end{equation}
\end{corollary}

  \begin{proof}
    The first inequality follows from the triangle inequality since
    \[
|\sqrt{2}e^{-it}-i|\ge \sqrt{2}-1.
\]
The second follows from
\[
  \begin{split}
    |\sqrt{2}e^{-it}-i|^2&=2+i\sqrt{2}e^{-it}-i\sqrt{2}e^{it}+1\\
    &=3+2\sqrt{2}\sin t\\
    &\le 3+2\sqrt{2}\\
    &=(\sqrt{2}+1)^2.
\end{split}
  \]
The last inequality follows from the Cauchy-Schwarz inequality.
\end{proof}

\subsection{A moment problem and extension problems}
From \eqref{marseille} we get
\begin{eqnarray}
  \nonumber
  f(m,0)&=&\frac{1}{\pi}\int_0^{2\pi}\left(\sum_{k=0}^m{m\choose k}\sqrt{2^k}e^{-ikt}(-i)^{m-k}\right)dM(t)\\
 \label{fmmm}
&=&2\sum_{k=0}^m\sqrt{2^k}{m\choose k}(-i)^{m-k}M_k,\quad m=0,1,\ldots
  \end{eqnarray}
where $M_k=\frac{1}{2\pi}\int_0^{2\pi}e^{-ikt}dM(t)$, $k=0,1,\ldots$ denote the moments of the spectral measure.

\begin{proposition} It holds that
  \begin{equation}
    \label{Mnfm}
    M_m=\frac{1}{\sqrt{2^m}}\sum_{k=0}^m{m\choose k}i^{m-k}f(k,0),\quad m=0,1,\ldots
    \end{equation}
\end{proposition}

\begin{proof} We have
  \[
\begin{split}
  M_m&=\frac{1}{2\pi}\int_0^{2\pi}e^{-imt}dM(t)\\
  &=\frac{1}{2\pi}\int_0^{2\pi}\left(\frac{\sqrt{2}e^{-it}-i+i}{\sqrt{2}}\right)^mdM(t)\\
  &=\frac{1}{2\pi}\frac{1}{\sqrt{2^m}}\sum_{k=0}^m{m\choose k}\int_0^{2\pi}(\sqrt{2}e^{-it}-i)^ki^{m-k}dM(t)\\
  &=\sum_{k=0}^ma_{m,k}f(k,0)
\end{split}
\]
with
\begin{equation}
  a_{m,k}=\frac{i^{m-k}{m\choose k}}{2\sqrt{2^m}},\quad k=0,\ldots, m.
\end{equation}
\end{proof}

The above relations suggest connections with the classical trigonometric moment problem.

\begin{proposition}
The main minors $F_N$ associated to $f(m,n)$ are non-negative if and only if the same holds for the Toeplitz matrices $T_N$
\begin{equation}
  \label{toeplitz98}
   2 T_N=L_NF_NL_N^*.
    \end{equation}
\end{proposition}
\begin{proof}
  Write as in the proof of the previous proposition
  \[
  e^{-it}=\frac{\sqrt{2}e^{-it}-i+i}{\sqrt{2}}.
    \]
    Then,
\begin{equation}
    \label{place-saint-michel}
  \begin{pmatrix}
    1\\ e^{-it}\\ \vdots\\ e^{-iNt}\end{pmatrix}=L_N\begin{pmatrix}1\\ \sqrt{2}e^{-it}-1\\ \vdots \\ (\sqrt{2}e^{-it}-1)^N\end{pmatrix}
\end{equation}
where $L_N$ is the invertible lower triangular matrix defined by
\[
(L_N)_{km}=\frac{1}{\sqrt{2^m}}{m\choose k} i^{m-k},\quad k=0,\ldots, m, \quad m=0,\ldots, N
  \]
From \eqref{place-saint-michel} we get
    \[
 \begin{pmatrix}   1\\ e^{-it}\\ \vdots\\ e^{-iNt}\end{pmatrix} \begin{pmatrix}   1\\ e^{-it}\\ \vdots\\ e^{-iNt}\end{pmatrix}^*
 =L_N\begin{pmatrix}1\\ \sqrt{2}e^{-it}-1\\ \vdots \\ (\sqrt{2}e^{-it}-1)^N\end{pmatrix}\begin{pmatrix}1\\ \sqrt{2}e^{-it}-1\\ \vdots \\ (\sqrt{2}e^{-it}-1)^N\end{pmatrix}^*L_N^*
\]
and integrating with respect to $dM(t)$ on $[0,2\pi]$ we get \eqref{toeplitz98}.
    \end{proof}

  \begin{corollary}
In the positive case the generating function converges in a neighborhood of the origin.
  \end{corollary}

  \begin{proof} We focus on the scalar case first.
    Let $(f(m,n))$ be a positive discrete analytic function, with a possibly divergent generating function, and define $M_0,M_2,\ldots$ via \eqref{Mnfm}. The corresponding Toeplitz matrices
    are non negative and hence $|M_m|\le M_0$ and the series $\varphi(z)= M_0+2\sum_{m=1}^\infty M_m z^m$ converges in the open unit disk. Thus the function $\Phi(z)=\varphi(\sigma^{-1}(z))$
    converges in $B(-i,\sqrt{2})$, and the result follows.\smallskip

    The $\mathbb C^{p\times p}$-valued case is treated by first considering the functions $c^*f(m,n)c$ for $c\in\mathbb C^p$ and using the polarization identity.
    \end{proof}

Recall that matrices with discrete analytic structure were introduced in Definition \ref{dafm}.

\begin{corollary}
 Equations \eqref{fmmm} and \eqref{Mnfm} give a one-to-one correspondence between positive matrices with discrete analytic structure
  and positive Toeplitz matrices.
  \end{corollary}

\begin{problem}
  Given a positive matrix with discrete analytic structure $F_N=(f(m,n))_{n,m=0}^N$ describe the set of all numbers ${\lambda}$ such that there is a positive discrete analytic matrix of the form
  \[
  F_{N+1}({\lambda})=  \begin{pmatrix}&&&\overline{{\lambda}}\\
      & F_N& &*\\
& & &\vdots\\
      {\lambda}&*&\cdots &*\end{pmatrix}
  \]
  where the $*$ denote entries to be computed from $F_N$ and ${\lambda}$. For a given ${\lambda}$, these entries are uniquely determined by the discrete Cauchy-Riemann equations, i.e. by
  \eqref{def-daf}.
\end{problem}

We state the next result without a detailed proof.

\begin{theorem}
  There is a one-to-one correspondence between the discrete analytic extension of a given
  positive matrix with discrete analytic structure and the corresponding positive extensions of the associated Toeplitz matrix.
\end{theorem}

\begin{proof}
  This follows from \eqref{toeplitz98}.
\end{proof}

The following subsection expands on these connections. The case of rational spectral function is considered in Section \ref{gohgoh}.

\subsection{Schur analysis and discrete analytic functions in $\Lambda_{++}$}
To every Carath\'eodory function is associated a discrete analytic function in $\Lambda_{++}$ of a special type (positivity of the principal matrices $F_N$) and a family of methods and problems
(called Schur analysis, or $J$-theory; see \cite{Dym_CBMS}) pertaining to
\begin{itemize}
\item Interpolation.

  \item Stationary stochastic processes.

  \item Inverse scattering and operator models.

    \item Schur coefficients.
    \end{itemize}
    We postpone the study of the connections between these notions of Schur analysis and the associated discrete analytic function
and contend ourselves to the study of reproducing kernel Hilbert spaces with reproducing kernel of the form
    \[
k_f(\lambda,\nu)=\frac{\Phi(\lambda)+\Phi(\nu)^*}{1+i\lambda-i\overline{\nu}-\lambda\overline{\nu}}.
    \]

    \begin{proposition}
      The reproducing kernel Hilbert space $\mathcal L(k_f)$ with reproducing kernel $k_f$ consists
      of functions of the form
\begin{equation}
  \mathsf F ({\lambda})=\frac{F(\sigma({\lambda}))}{\pi(1+i{\lambda})}\, ,\, F\in\mathcal L(\varphi)
\end{equation}
with norm $\|\mathsf F\|=\|F\|_{\mathcal L(\varphi)}$, where $\mathcal L(\varphi)$ denotes the reproducing kernel Hilbert space with reproducing kernel
$\frac{\varphi(\lambda)+(\varphi(\nu))^*}{1-\lambda\overline{\nu}}$.
      \end{proposition}

      \begin{proof} The claim follows from
        \begin{equation}
          \begin{split}
  k_f({\lambda},{\nu})&=\frac{1}{\pi(1+i{\lambda})(1-i\overline{\nu})}\int_0^{2\pi}\frac{dM(t)}{(1-\sigma({\lambda})e^{-it})(1-\overline{\sigma({\nu})}e^{it})}\\
  &=\frac{1}{\pi}\int_0^{2\pi}\frac{dM(t)}{(a({\lambda})-b({\lambda})e^{-it})(\overline{a({\nu})}-\overline{b({\nu})}e^{it})}.
  \end{split}
\end{equation}
\end{proof}

It is of interest to characterize intrinsically reproducing kernel spaces Hilbert spaces with reproducing kernel $\frac{\Phi(\lambda)+\Phi(\nu)^*}{1+i\lambda-i\overline{\nu}-\lambda\overline{\nu}}.$
We give below a general characterization for arbitrary denominators of the form $a(\lambda)\overline{a(\nu)}-b(\lambda)\overline{b(\nu)}$, where the functions $a$ and $b$ are analytic in a connected
open set $\Omega$, and such that the corresponding sets $\Omega_+$ and $\Omega_-$ defined in \eqref{1-20-20-20} are non-empty (then $\Omega_0$ is also non-empty).
With $\sigma=b/a$ one defines (see \cite[(3.4), (3.5), p. 8]{ad-jfa})
\begin{align}
\label{4-17}
  (r(a,b,\alpha))f({\lambda})&=\frac{a({\lambda})f({\lambda})-a(\alpha)f(\alpha)}{a(\alpha)b({\lambda})-b(\alpha)a({\lambda})}\\
  (r(b,a,\alpha))f({\lambda})&=\frac{b({\lambda})f({\lambda})-b(\alpha)f(\alpha)}{b(\alpha)a({\lambda})-a(\alpha)b({\lambda})}
\label{4-18}
\end{align}
  Note  that

  \begin{equation}
    a(\alpha)r(a,b,\alpha)-b(\alpha)r(b,a,\alpha)=Id\quad (\mbox{\text{\rm the identity operator}})
\end{equation}

For fixed $t\in\mathbb R$ it holds that (see equation (4.2) in \cite{ad-jfa}):
  \begin{equation}
    \left(r(a,b,\alpha)\frac{1}{a-be^{-it}}\right)({\lambda})=\frac{e^{-it}}{a(\alpha)-b(\alpha)e^{-it}}
    \frac{1}{a({\lambda})-b({\lambda})e^{-it}}
  \end{equation}
  and
  \begin{equation}
    \left(r(b,a,\alpha)\frac{1}{a-be^{-it}}\right)({\lambda})=\frac{-1}{a(\alpha)-b(\alpha)e^{-it}}
    \frac{1}{a({\lambda})-b({\lambda})e^{-it}}
  \end{equation}
  When $a(\lambda)=1$ and $b(\lambda)=1$, we have $r(a,b,\alpha)=R_\alpha$, with
  \begin{equation}
    (R_\alpha f)(\lambda)=\begin{cases}\,\dfrac{f(\lambda)-f(\alpha)}{\lambda-\alpha}\,\,\lambda\not=\alpha\\
    f^\prime(\alpha),\,\,\,\hspace{11mm} \lambda=\alpha.
\end{cases}
\end{equation}

\begin{definition}
We say that $\Omega\subset \mathbb C$ is symmetric with respect to $C(-i,\sqrt{2})$ if the following condition holds: $\lambda\in\Omega$ if and only if $\lambda^\prime\in\Omega$ where $\lambda^\prime=\frac{1-i\overline{\lambda}}{\overline{\lambda}-i}$, i.e. (see Lemma \ref{lemmalrlr}) if and only if $\lambda$ and $\lambda^\prime$ are related by $\sigma(\lambda)\overline{\sigma(\lambda^\prime)}=1$.
\end{definition}
\begin{theorem}
  Let $\mathfrak H$ be a reproducing kernel Hilbert space of $\mathbb C^p$-valued functions analytic in an open set symmetric with respect to $C(-i,\sqrt{2})$ and such that
  \begin{equation}
    a(\alpha)b^\prime(\alpha)-b(\alpha)a^\prime(\alpha)\not=0,\quad \alpha\in\Omega.
    \end{equation}
Then the reproducing kernel of $\mathfrak H$ is of the form
\[
  \frac{\Phi(\lambda)+\Phi(\nu)^*}{1+i\lambda-i\overline{\nu}-\lambda\overline{\nu}}
\]
for some $\mathbb C^{p\times p}$valued function analytic in $\Omega$ if and only if it is invariant under both $r(b,a,\alpha)$ and $r(a,b,\alpha)$ for all $\alpha\in\Omega$ and it holds that
  \begin{equation}
    \label{structure}
    [r(b,a,\alpha)\mathsf F, r(b,a,\beta)\mathsf G]-
     [r(a,b,\alpha)\mathsf F, r(a,b,\beta)\mathsf G]=0
   \end{equation}
   for all $\mathsf F,\mathsf G\in\mathfrak H$.
\end{theorem}

\begin{proof}
  We will prove this result as a special case of \cite[Theorem 4.1, p. 12]{ad-jfa} and proceed in a number of steps.\\

  STEP 1: {\sl We take $J=\begin{pmatrix}0&I_p\\I_p&0\end{pmatrix}$ and define
      \[
      \mathfrak H_e=\left\{\mathsf F_e(\lambda)=\begin{pmatrix}\mathsf F(\lambda)\\ 0_{n\times 1}\end{pmatrix}\,;\,\mathsf F\in\mathfrak H\,\,{\rm and}\,\,
      \|\mathsf F_e\|_{\mathfrak H_e}=\|\mathsf F\|_{\mathfrak H}
          \right\}.
          \]
          Then, the structural identity \eqref{structure} can be rewritten as
\begin{equation}
    \label{structure-1234}
    [r(b,a,\alpha)\mathsf F_e, r(b,a,\beta)\mathsf G_e]-
     [r(a,b,\alpha)\mathsf F_e, r(a,b,\beta)\mathsf G_e]=\mathsf G_e(\beta)^*J\mathsf F_e(\alpha).
   \end{equation}
 }

This is because the right hand side of \eqref{structure-1234} is equal to $0$ since the range of the functions in $\mathfrak H_e$ is $J$-neutral with the above choice of $J$.\\

STEP 2: {\sl There is a matrix-valued function $\Theta$ analytic in $\Omega$ and such that such that the reproducing kernel $K(\lambda, \nu)$ of $\mathfrak H_e$ is
\[
K(\lambda,\nu)=\frac{J-\Theta(\lambda)J\Theta(\nu)^*}{a(\lambda)\overline{a(\nu)}-b(\lambda)\overline{b(\nu)}}.
\]}

This follows by applying \cite[Theorem 4.3 p. 16]{ad-jfa}. \smallskip

The statement in \cite{ad-jfa} requires that, for some point $\mu$ having a symmetric point $\mu^\prime\in\Omega$.
one can write
\begin{equation}
  \label{eje}
J-\rho_\mu(\mu)K(\mu,\mu)=E^*JE
  \end{equation}
  for some $p\times p$-invertible matrix $E$. This condition is automatically satisfied in the Hilbert space case,
  as has been shown in the disk case by Rovnyak in \cite{rov-66} and Ball in \cite[Lemma 2 p. 251]{ball-contrac} in the line case.
  A proof for the case of a general pair of functions $(a,b)$ is given in
  \cite[Theorem 4.3 p. 16]{ad-jfa}. The crux of the proof of the present theorem is that $E$ has a specific form here, as explained in the next step.\smallskip

STEP 3: {\sl The matrix $E$ in \eqref{eje} can be taken of the form $\begin{pmatrix}V&I_p\\I_p&0\end{pmatrix}$ for some self-adjoint matrix $V\in\mathbb C^{p\times p}$.}\smallskip

Applying to the reproducing kernel of $\mathfrak H_e$ the Zaremba-Bergman reproducing kernel formula\footnote{
(the paper \cite{franek_hbot1} sets the record straight on the history of the
formula. We choose to keep the name Bergman)}
(see \cite{Zaremba,franek_hbot1,Bergman22}) in terms of an orthonormal basis, we get
\begin{equation}
  \label{rtyuio}
\frac{J-\Theta(\lambda)J\Theta(\nu)^*}{a(\lambda)\overline{a(\nu)}-b(\lambda)\overline{b(\nu)}}=\sum_{n\in A}
\begin{pmatrix} e_n(\lambda)&0\end{pmatrix}\begin{pmatrix}e_n(\nu)^*\\ 0\end{pmatrix}=\begin{pmatrix}k(\lambda,\nu)&0\\0&0\end{pmatrix},
\end{equation}
where $\left\{\begin{pmatrix}e_n\\ 0\end{pmatrix},\, n\in A\right\}$ is {\sl any} orthonormal basis of
$\mathfrak H_e$. The analyticity of the kernel
  in $\lambda$ and $\overline{\nu}$ implies that the associated reproducing kernel Hilbert space is separable, and so we can take
  $A=\mathbb N$ as index set for the orthonormal basis, or $1,\ldots, N$ when the space has finite dimension $N$. Equation \eqref{rtyuio}
  implies that only the $p\times p$ left upper block of $K(\lambda,\nu)$ is not vanishing.
  Equation \eqref{eje} can thus be written as
  \begin{equation}
    \begin{pmatrix}U&I_p\\I_p&0\end{pmatrix}=E^*JE
  \end{equation}
  with $U=-\rho_\mu(\mu)\sum_{n\in A}e_n(\mu)e_n(\mu)^*$. It suffices then to take $E=\begin{pmatrix}U/2&I_p\\I_p&0\end{pmatrix}$.
      One has $E=E^*$ and
      \[
E^*JE=\begin{pmatrix}U/2&I_p\\I_p&0\end{pmatrix}J\begin{pmatrix}U/2&I_p\\I_p&0\end{pmatrix}=\begin{pmatrix}U&I_p\\I_p&0\end{pmatrix}.
        \]

STEP 4: {\sl   $\Theta$ is of the form
\[
\begin{pmatrix}\Phi(\lambda)&I_p\\
  I_p&0\end{pmatrix}
  \]
  and so
  \[
  \frac{J-\Theta(\lambda)J\Theta(\nu)^*}{a(\lambda)\overline{a(\nu)}-b(\lambda)\overline{b(\nu)}}=
  \begin{pmatrix}\frac{\Phi(\lambda)+\Phi(\nu)^*}{a(\lambda)\overline{a(\nu)}-b(\lambda)\overline{b(\nu)}}&0\\0&0\end{pmatrix}
    \]
    and the result follows since the reproducing kernel of $\mathfrak H$ is
    \[
\begin{pmatrix}I_p&0\end{pmatrix}\frac{J-\Theta(\lambda)J\Theta(\nu)^*}{a(\lambda)\overline{a(\nu)}-b(\lambda)\overline{b(\nu)}} \begin{pmatrix}I_p\\0\end{pmatrix}.
    \]}

  To check the above step it is enough to use \cite[(4.4), p. 12]{ad-jfa} with $E$ as in Step 3. More precisely, and with $k(\lambda,\mu)$ as in \eqref{rtyuio}:
  \[
    \begin{split}
      \Theta(\lambda)&=\left(J-\rho_{\mu^\prime}(\lambda)K(\lambda,\mu^\prime)\right)JE^*\\
      &=\begin{pmatrix}0&I_p\\ I_p&0\end{pmatrix}
      -\begin{pmatrix}\rho_{\mu^\prime}(\lambda)k(\lambda,\mu^\prime)
        &0\\0&0\end{pmatrix}\begin{pmatrix}0&I_p\\ I_p&0\end{pmatrix}\begin{pmatrix}U/2&I_p\\I_p&0\end{pmatrix}\\
      &=\begin{pmatrix}-\rho_{\mu^\prime}(\lambda)k(\lambda,\mu^\prime)&I_p\\I_p&0\end{pmatrix},
      \end{split}
    \]
    and this concludes the proof.
\end{proof}

We note that for the disk case we have
\[
[(I+\alpha R_\alpha)F,(I+\beta R_\beta)G]-[R_\alpha F,R_\beta G]=0\quad (\mbox{or, more generally, $G(\beta)^*JF(\alpha)$})
\]
i.e.
\begin{equation}
  [F,G]+\alpha[R_\alpha F,G]+\overline{\beta}[F,R_\beta G]-(1-\alpha\overline{\beta})[R_\alpha F,R_\beta G]=0 \quad (\mbox{or, more generally, $G(\beta)^*JF(\alpha)$}).
  \end{equation}
Here we have:

\begin{proposition} In the case of $a(\lambda)=1+i\lambda$ and $b(\lambda)=\sqrt{2}\lambda$
  we have
  \begin{align}
    (r(a,b,\alpha)f)(\lambda)&=\frac{i}{\sqrt{2}}f(\lambda)+\left(\frac{1+i\alpha}{\sqrt{2}}\right)(R_\alpha f)(\lambda)\\
    (r(b,a,\alpha)f)(\lambda)&=-f(\lambda)-\alpha R_\alpha f(\lambda)
                               \end{align}
  and the structural identity is equal to
  \begin{equation}
    \label{4-29}
    \begin{split}
      \hspace{-1cm}[F,G]+(-i+\overline{\beta})[F,R_\beta G]+(i+\alpha)[R_\alpha F,G]-(1-i\overline{\beta}+i\alpha-\alpha\overline{\beta})[R_\alpha F,R_\beta G]=\\
      &\hspace{-6cm}=\begin{cases}\,\,0\,\, \mbox{\hspace{18mm}for $\mathcal L(\Phi)$ spaces}\\
G(\beta)^*JF(\alpha)\,\, \mbox{for $\mathcal H(\Theta)$ spaces.}\end{cases}
\end{split}
          \end{equation}
\end{proposition}

\begin{proof} We first compute the operators $r(a,b,\alpha)$ and $r(b,a,\alpha)$ using \eqref{4-17} and \eqref{4-18} and get
  \[
    \begin{split}
      (r(a,b,\alpha)f)(\lambda)&=\frac{(1+i\lambda)f(\lambda)-(1+i\alpha)f(\alpha)}{(1+i\alpha)\sqrt{2}\lambda-\sqrt{2}\alpha(1+i\lambda)}\\
      &=\frac{f(\lambda)-f(\alpha)+i(\lambda f(\lambda)-\alpha f(\alpha))}{\sqrt{2}(\lambda-\alpha)}\\
            &=\frac{f(\lambda)-f(\alpha)+i(\lambda -\alpha+\alpha)f(\lambda)-i\alpha f(\alpha)}{\sqrt{2}(\lambda-\alpha)}\\
            &=\frac{i}{\sqrt{2}}f(\lambda)+\left(\frac{1+i\alpha}{\sqrt{2}}\right)(R_\alpha f)(\lambda)
          \end{split}
        \]
and

        \[
    \begin{split}
      (r(b,a,\alpha)f)(\lambda)&=\frac{\sqrt{2}\lambda f(\lambda)-\sqrt{2}\alpha f(\alpha)}{\sqrt{2}\alpha(1+i\lambda)-(1+i\alpha)   \sqrt{2} \lambda}\\
      &=\frac{\lambda f(\lambda)-\alpha f(\alpha)}{\alpha- \lambda}\\
    &=\frac{(\lambda-\alpha+\alpha) f(\lambda)-\alpha f(\alpha)}{\alpha- \lambda}\\
&=      -f(\lambda)-\alpha R_\alpha f(\lambda).
\end{split}
\]
So the first two formulas in the statement follow. Then \eqref{structure} becomes:
\[
  \begin{split}
    \left(R_\beta^*\left(\frac{1-i\overline{\beta}}{\sqrt{2}}-\frac{i}{\sqrt{2}}I\right)\right)\left(\left(\frac{1+i\overline{\alpha}}{\sqrt{2}}+\frac{i}{\sqrt{2}}I\right)R_\alpha\right)-
    (\overline{\beta}R_\beta^*+I)(I+\alpha
    R_\alpha)&=\\
    &\hspace{-12cm}=R_\beta^*\left(\left(\frac{1-i\overline{\beta}}{\sqrt{2}}\right)
      \left(\frac{1+i\overline{\alpha}}{\sqrt{2}}\right)-\overline{\beta}\alpha\right)R_\alpha-R_\beta^*\left(\overline{\beta}-\frac{i}{\sqrt{2}}\frac{1-i\overline{\beta}}{\sqrt2}\right)-\\
    &\hspace{-11.5cm}-\left(\alpha+\frac{-i}{\sqrt{2}}\frac{1+i\alpha}{\sqrt{2}}\right)R_\alpha+\frac{I}{2}-I\\
    &\hspace{-12cm}=R_\beta^*R_\alpha\frac{1-i\overline{\beta}+i\alpha-\alpha\overline{\beta}}{2}+R_\beta^*\frac{i-\overline{\beta}}{2}+R_\alpha\frac{-i-\alpha}{2}-\frac{I}{2}
    \end{split}
\]
and the assertion follows.
  \end{proof}

  \begin{remark}{\rm
When we take $\alpha=\beta=0$ in \eqref{4-29} we obtain
\begin{equation}
  \label{4-28-000}
I+i(R_0-R_0^*)-R_0^*R_0=C^*JC
  \end{equation}

to compare with
\[
I-R_0^*R_0=C^*JC
\]
for the disk case, and
\[
iR_0-iR_0^*=C^*JC
\]
for the half-plane case.}
\end{remark}

\begin{remark}
  {\rm
    Positive discrete analytic functions in $\Lambda_{++}$ are in a one-to-one correspondence with functions analytic and contractive in the open unit disk (at the exception of the function
    $s(\lambda)\equiv -1$) via the formula
  \[
\frac{1-\Phi(\sigma^{-1}(\lambda))}{1+\Phi(\sigma^{-1}(\lambda))}.
    \]
}
\end{remark}

The discrete analytic function which  corresponds to a degree one Blaschke is computed in Example \ref{blabla}.

  \section{The rational case}
  \setcounter{equation}{0}
  In the positive case, the spectral function $dM(t)$ is the weak-$*$ limit of the measures ${\rm Re}\, \varphi(re^{it})dt$ as $r\uparrow 1$. This makes use of the
  Banach-Alaoglu theorem, or of  Helly's theorem; see e.g. \cite[pp. 19-25]{MR48:904} for a discussion using Helly's theorem. It follows from the above
  characterization of $dM$ that $\varphi$ is rational and lossless if and only if
  $dM$ is jump measure with a finite number of jumps, and that, when $\varphi$ is rational and  has no singularities on the boundary, the measure $M$ is absolutely
  continuous with respect to measure Lebesgue and
  \begin{equation}
    {\rm Re}\, \varphi(e^{it})=M^\prime(e^{it})
    \label{re-phi}
\end{equation}
    In this section we obtain formulas for the values $f(m,0)$ both in the lossless
    rational case and in the case just described above by \eqref{re-phi}.
Extension to right-upper quarter-plane is done using unitary dilations theory.\smallskip

One could extend these examples to the case where the
rational function is not analytic at infinity, by using other kinds of realizations (see for instance \cite[(3.10) p. 398 and Theorem 8.2 p. 422]{gk1} for the latter).
As mentioned in Section \ref{sub-real}, this will not be done here.\\

\subsection{Rational spectral function}
\label{gohgoh}
We first compute $\varphi$ when $M$ is absolutely continuous with Lebesgue
measure. Let us consider the minimal realization
\[
M^\prime(\lambda)=D+C(\lambda I_{2\ell}-A)^{-1}B
\]
of $M^\prime$ where
$A,B,C,D$ given by \eqref{prod111}-\eqref{prod333}. In particular, $M^\prime$ is assumed invertible at the origin and analytic at infinity. Since

\[
\varphi(\lambda)=\frac{1}{2\pi}\int_0^{2\pi}\frac{e^{it}+\lambda}{e^{it}-\lambda}M^\prime(e^{it})dt
=\frac{1}{2\pi}\int_0^{2\pi}M^\prime(e^{it})dt+2\lambda+\frac{1}{2\pi}\int_0^{2\pi}\frac{1}{e^{it}-\lambda}M^\prime(e^{it})dt
\]
we need first to compute

  \[
    \begin{split}
      \frac{1}{2\pi} \int_0^{2\pi}\frac{1}{e^{it}-\lambda}(e^{it}I_{2\ell}-A)^{-1}dt&=\\
      &\hspace{-3cm}=\frac{1}{2\pi} \int_0^{2\pi}\frac{1}{e^{it}-\lambda}\begin{pmatrix}(e^{it}I_{\ell}-a)^{-1}dt
        &-(e^{it}I_{\ell}-a)^{-1}bb^*a^{-*}(e^{it}I_{\ell}-a^{-*})^{-1}\\
        0&(e^{it}I_{\ell}-a^{-*})^{-1}\end{pmatrix}dt\\
        &\hspace{-3cm}=\begin{pmatrix}\spadesuit_1&\spadesuit_2\\0&\spadesuit_3\end{pmatrix}.
  \end{split}
  \]
  The block entries $\spadesuit_j$, $j=1,2,3$ are computed as follows:

  \[
      \spadesuit_1=\frac{1}{2\pi i}\int_{|\zeta|=1}\frac{1}{(\zeta-z)\zeta}(\zeta {I_{\ell}}-a)^{-1}d\zeta=0
        \]
    by Abel’s relation on the sum of the residues of a rational function; see for instance \cite[p. x]{MR1704477} for the latter. Then we have:\\

    \[
      \begin{split}
        \spadesuit_3&=\frac{-1}{2\pi}\int_0^{2\pi}\frac{1}{e^{it}-\lambda}(1-a^*e^{it})^{-1}a^*dt\\
        &=\frac{-1}{2\pi}\int_0^{2\pi}\left(\sum_{m=0}^\infty e^{-(m+1)t}\lambda^m\right)\left(\sum_{n=0}^\infty e^{int}a^{*n}\right)a^*dt\\
        &=-a^{*2}(I_{\ell}-a^*\lambda)^{-1}.
      \end{split}
    \]
    To compute $\spadesuit_2$ we first note that
    \begin{equation}
      (e^{it}I_{\ell}-a)^{-1}\cdot(e^{it}I_{\ell}-\lambda I_{\ell})^{-1}=(a-\lambda I_{\ell})^{-1}\left((e^{it}I_{\ell}-a)^{-1}-(e^{it}I_{\ell}-\lambda I_{\ell})^{-1}\right).
      \label{resolv}
    \end{equation}
    Therefore, we can write:
    \[
      \begin{split}
        \spadesuit_2&=\frac{1}{2\pi}\int_0^{2\pi}(e^{it}I_{\ell}-a)^{-1}bb^*(I_{\ell}-e^{it}a^*)^{-1}\frac{1}{e^{it}-\lambda}dt\\
        &=\frac{1}{2\pi}\int_0^{2\pi}(a-\lambda I_{\ell})^{-1}(e^{it}I_{\ell}-a)^{-1}bb^*(I_{\ell}-e^{it}a^*)^{-1}dt+\\
        &\hspace{5mm}+\frac{1}{2\pi}\int_0^{2\pi}(a-\lambda I_{\ell})^{-1}bb^*(I_{\ell}-e^{it}a^*)^{-1}\frac{1}{e^{it}-\lambda}dt\\
        &=\spadesuit_{21}+\spadesuit_{22},
    \end{split}
  \]
where
    \[
    \begin{split}
      \spadesuit_{21}&=(a-\lambda I_{\ell})^{-1}\int_0^{2\pi}\left(\sum_{m=0}^\infty e^{-(m+1)t}a^m\right)bb^*\left(\sum_{n=0}^\infty e^{int}a^{*n}\right)dt\\
&=(a-\lambda I_{\ell})^{-1}Xa^*
    \end{split}
  \]
and
  \[
    \begin{split}
      \spadesuit_{22}&=-(a-\lambda I_{\ell})^{-1}\int_0^{2\pi}\left(\sum_{m=0}^\infty e^{-(m+1)t}\lambda^m\right)bb^*\left(\sum_{n=0}^\infty e^{int}a^{*n}\right)dt\\
      &=-(a-\lambda I_{\ell})^{-1} bb^*a^*(I_{\ell}-\lambda a^*)^{-1}.
    \end{split}
  \]

  Using the Stein equation \eqref{stein-equa}, we get
  \begin{equation}
    \begin{split}
      \spadesuit_2&=(a-\lambda I_{\ell})^{-1}\left(X-bb^*(I_{\ell}-a^*\lambda)^{-1}\right)a^*\\
      &=(a-\lambda I_{\ell})^{-1}\left(X-bb^*-Xa^*\lambda\right)a^*(I_{\ell}-\lambda a^*)^{-1}\\
      &=(a-\lambda I_{\ell})^{-1}\left(aXa^*-Xa^*\lambda\right)a^*(I_{\ell}-\lambda a^*)^{-1}\\
      &=Xa^{*2}(I_{\ell}-\lambda a^*)^{-1}.
      \end{split}
    \end{equation}

    Then we have
    \begin{equation}
      \frac{1}{2\pi} \int_0^{2\pi}\frac{1}{e^{it}-\lambda}(e^{it}I_{2\ell}-A)^{-1}dt=\begin{pmatrix}0& Xa^{*2}(I_{\ell}-\lambda a^*)^{-1}\\0& a^{*2}(I_{\ell}-\lambda a^*)^{-1}\end{pmatrix}
    \end{equation}
    and
    \[
\frac{1}{2\pi}\int_0^{2\pi}M^\prime(e^{it})dt=-\lim_{\lambda\rightarrow\infty} \frac{\lambda}{2\pi} \int_0^{2\pi}\frac{1}{e^{it}-\lambda}(e^{it}I_{2\ell}-A)^{-1}dt=\begin{pmatrix}0& Xa^*\\0& -a^*\end{pmatrix}.
    \]
    So,
    \[
    \frac{1}{2\pi}\int_0^{2\pi}M^\prime(e^{it})dt+
    2\lambda\frac{1}{2\pi} \int_0^{2\pi}\frac{1}{e^{it}-\lambda}(e^{it}I_{2\ell}-A)^{-1}dt=
    \begin{pmatrix}0& Xa^*(I_{\ell}+\lambda a^*)(I_{\ell}-\lambda a^*)^{-1}\\0& -a^*(I_{\ell}+\lambda a^{*})(I_{\ell}-\lambda a^*)^{-1}\end{pmatrix}
    \]

    and, with $B,C,D$ given by \eqref{prod222}-\eqref{prod333}
    \begin{equation}
      \label{no-garentee}
      \begin{split}
        \varphi(\lambda)&=D+C
\begin{pmatrix}0& Xa^*(I_{\ell}+\lambda a^*)(I_{\ell}-\lambda a^*)^{-1}\\0& -a^*(I_{\ell}+\lambda a^{*})(I_{\ell}-\lambda a^*)^{-1}\end{pmatrix}
        B\\
        &=d(d^*-b^*a^{-*}c^*)+(cX+db^*a^{-*})(I_{\ell}+\lambda a^*)(I_{\ell}-\lambda a^*)^{-1}c^*,
      \end{split}
      \end{equation}
      which is \eqref{dnot0}. When $d=0$ we get back to \eqref{form02222022a}, namely
    \begin{equation}
      \label{ca-marche}
    \varphi(\lambda)=cX(I_\ell+\lambda a^*)(I_\ell-\lambda a^*)^{-1}c^*=cX(a^{-*}+\lambda I_{\ell})(a^{-*}-
    \lambda I_\ell)^{-1}c^* .
    \end{equation}

\subsection{The values $f(m,0)$}
We assume now that the spectral measure is absolutely continuous with respect to the Lebesgue measure, with a  rational spectral density, regular
on the unit circle and at infinity. In the statement, note that $A,B,C$ and $D$ are not arbitrary, since $R(e^{it})>0$ and $R$ is assumed
regular  on the unit circle. In view of the assumed minimality the realization of $R$
is similar to a realization of the form \eqref{prod111}-\eqref{prod333}. Furthermore, as earlier in the paper, we write
(with an abuse of notation) $((I_n-P)A)^{-1}$  for $((I_N-P)A|_{{\rm ran}\,I_N-P})^{-1}(I_N-P)$.

\begin{theorem}
  Let the spectral density be $\mathbb C^{p\times p}$-valued, rational, with minimal realization
  $M^\prime(e^{it})=R(e^{it})dt$ with $R(\lambda)=D+C(\lambda I_N-A)^{-1}B$.
  Then
  \begin{equation}
    \label{inverse4}
f(m,0)=D-C(\sqrt{2}((I_N-P)A)^{-1}-iI_N)^mB,\quad m=0,1,\ldots
    \end{equation}
\end{theorem}
  \begin{proof}
We use the representation \eqref{coeff-fourier} for the Fourier coefficients of the weight and get:
\[
  \begin{split}
    f(m,0)&=
    \frac{1}{2\pi}\int_0^{2\pi}(\sqrt{2}e^{-it}-i)^m\left(\sum_{u\in\mathbb Z} R_ue^{iut}\right)dt\\
    &=\frac{1}{2\pi}\int_0^{2\pi}
  \left(\sum_{k=0}^m2^{k/2}e^{-ikt}(-i)^{m-k}{m\choose k}\right)\left(\sum_{u\in\mathbb Z} R_ue^{iut}\right)dt\\
  &=\sum_{k=0}^m 2^{k/2}(-i)^{m-k}{m\choose k}R_{k}\\
  &=(D-C((I_N-P)A)^{-1}(I_N-P)B)(-i)^m-\\
  &\hspace{5mm}-\sum_{k=1}^m C((I_N-P)A)^{-k-1}(I_N-P)B 2^{k/2}(-i)^{m-k} .
  \end{split}
\]

Since
\[
\sum_{k=0}^m ((I_N-P)A)^{-k}2^{k/2}(-i)^{-k}{m\choose k}=(\sqrt{2}((I_N-P)A)^{-1}(-i)^{-1}+I_N)^m,
\]
we have
\[
  \begin{split}
    f(m,0)&=D-(-i)^mC\left(\sum_{k=0}^m ((I_N-P)A)^{-k}2^{k/2}(-i)^{-k}{m\choose u}\right)C\\
    &=D-(-i)^mC(\sqrt{2}((I_N-P)A)^{-1}(-i)^{-1}+I_N)^mB\\
    &=D-C(\sqrt{2}((I_N-P)A)^{-1}-iI_N)^mB,
\end{split}
\]
and the statement follows.
\end{proof}

  \begin{remark}{\rm One could have used the functional calculus associated to $A$ to get to the same result.}
  \end{remark}

Using the representations \eqref{r0}-\eqref{r-k} we get:

\begin{proposition}
  In the notation of the previous theorem, let $d+c(zI_{\ell}-a)^{-1}b$ be a minimal realization of the left spectral factor of $R(\lambda)$.
  Then:
  \begin{equation}
f(m,0)=(-i)^md(d^*-b^*a^{-*}c^*)+(db^*a^{-*}+cX)(\sqrt{2}a^*-iI_{\ell})^mc^*,\quad m=0,1,\ldots
    \end{equation}
and in particular, for $d=0$
\begin{equation}
  \label{d=0}
f(m,0)=cX(\sqrt{2}a^*-iI_{\ell})^mc^*,\quad m=0,1,\ldots
\end{equation}
and for $d=ca^{-1}b$,
  \begin{equation}
\label{le-louvre}
    f(m,0)=c(Xa^{-*}+a^{-1}X)(\sqrt{2}a^*-iI_{\ell})^mc^*,\quad m=0,1,\ldots
  \end{equation}
  \label{prop5-3}
\end{proposition}

\begin{proof} 
The proof is based on the following computations:
  \[
    \begin{split}
      \sum_{k=0}^m2^{k/2}(-i)^{m-k}{m\choose k}R_k&=2^{m/2}(-i)^m{m\choose 0}(dd^*+cXc^*)+\\
      &\hspace{5mm}+\sum_{k=1}^m2^{k/2}(-i)^{m-k}{m\choose k}(db^*+cXc^*)a^{*(k-1)}c^*\\
      &=(-i)^m{m\choose 0}(dd^*+cXc^*)-(-i)^m{m\choose 0}(db^*+cXa^*)a^{-*}c^*+\\
      &\hspace{5mm}+\sum_{k=0}^m2^{k/2}(-i)^{m-k}{m\choose k}(db^*+cXc^*)a^{*(k-1)}c^*\\
      &=(-i)^md(d^*-b^*a^{-*}c^*)+(db^*+cXa^*)(\sqrt{2}a^*-iI_{\ell})^ma^{-*}c^*.
      \end{split}
    \]
  \end{proof}
  Here too, we can define inverse problems using \eqref{inverse4}, consisting at recovering $A(I_N-P)$, $B$, $C$ and $D$, or $a,b,c$ and $d$, from $f(0,0),\ldots, f(N,0)$. This is the
  partial realization problem. See \cite{gkl-scl} for applications.\smallskip

  When $d=0$ we wish to extend \eqref{d=0} to $\Lambda_{++}$. In general $a^*\not=a$, but $a$ is a contraction. So it has a unitary dilation $U_a$ to a larger Hilbert space $\mathfrak H$,
  meaning that
  \begin{align}
    PU_a^n|_{\mathbb C^n}&=a^n,\quad n=0,1,\ldots\\
      PU_a^{-n}|_{\mathbb C^n}&=a^{*n},\quad n=-1,-2,\ldots
        \end{align}
  and we have as an extension
  \begin{equation}
    f(m,n)=CXP_{\mathbb C^N}(\sqrt{2}U_a^*-iI_{\mathfrak H})^m(\sqrt{2}U_a+iI_{\mathfrak H})^n|_{\mathbb C^N}C^*
  \end{equation}

  The above extension is not symmetric. To get a symmetric extension we rewrite \eqref{d=0} as

    \begin{equation}
\label{le-louvre-2}
f(m,0)=cX^{1/2}(\sqrt{2}X^{-1/2}a^*X^{-1/2}-iI_{\ell})^mX^{1/2}c^*,\quad m=0,1,\ldots,
\end{equation}
and apply the unitary dilation theorem to ${X}^{-1/2}a^*{X}^{-1/2}$.
\subsection{Rational lossless case}
By rational lossless case we mean the case where the $\mathbb C^{p\times p}$-valued function $\varphi$ in \eqref{kphi} is rational and satisfies
\[
\varphi(\lambda)+\varphi(\lambda)^*=0,\quad |\lambda|=1\quad(\mbox{\text{at the possible exception of poles on the unit circle}}),
\]
and
\begin{equation}
  \label{enst}
{\rm Re}\, \varphi(\lambda)\ge 0, \quad |\lambda|<1.
\end{equation}
Let $\varphi(\lambda)=D+C(\lambda I_N-A)^{-1}B$ be a minimal realization of $\varphi$, characterized by \cite[Theorem 5.1 p. 215]{ag}
(see Theorem \ref{th-ag-p-215} above). We take without loss of generality $H=I_N$ since $H>0$ in view of \eqref{enst} and \eqref{comp1}, and \eqref{comp22} becomes
\begin{equation}
  \label{phi-strasbourg}
\varphi({\lambda})=\frac{1}{2}C(A-{\lambda}I_N)^{-1}(A+{\lambda}I_N)C^*+iX,\quad X=X^*\quad{\rm and}\quad AA^*=A^*A=I_N.
\end{equation}
We set $X=0$, i.e. $\varphi(0)\ge 0$.

\begin{proposition}
  In the above notation,
      \begin{equation}
      f(m,0)=C(\sqrt{2}A^{-1}-iI_N)^mC^*,\quad m=0,1,\ldots
      \label{newform54321}
      \end{equation}

\end{proposition}

\begin{proof}

The same computations as above give:
\[
  \begin{split}
    \Phi({\lambda})&=\varphi(\sigma({\lambda}))\\
    &=\frac{1}{2}C\left(A-\frac{\sqrt{2}{{\lambda}}}{1+i{\lambda}}I\right)^{-1}\left(A+\frac{\sqrt{2}{{\lambda}}}{1+i{\lambda}}I\right)C^*\\
        &=\frac{1}{2}C\left((1+i{\lambda})A-\sqrt{2}{{\lambda}}I\right)^{-1}\left((1+i{\lambda})A+\sqrt{2}{{\lambda}}I\right)C^*\\
        &=\frac{1}{2}C(A+{\lambda}(iA-\sqrt{2})I)^{-1}(A+{\lambda}(iA+\sqrt{2})I)C^*\\
        &=\frac{1}{2}CC^*+\frac{1}{2}C\left((A+{\lambda}(iA-\sqrt{2}I))^{-1}(A+{\lambda}(iA+\sqrt{2})I)-I\right)C^*\\
        &=\frac{1}{2}CC^*+\sqrt{2}{\lambda}C(A+{\lambda}(iA-\sqrt{2}I))^{-1}C^*\\
        &=\frac{1}{2}CC^*+\sqrt{2}{\lambda}CA^{-1}(I-{\lambda}(\sqrt{2}A^{-1}-iI))^{-1}C^*\\
        &=\frac{1}{2}CC^*+\sum_{n=0}^\infty {\lambda}^{n+1}\sqrt{2}CA^{-1}(\sqrt{2}A^{-1}-iI)^nC^*\\
        &=\sum_{n=0}^\infty\Phi_n{\lambda}^n.
    \end{split}
  \]

  By \eqref{nantes}

  \[
    \begin{split}
      f(m,0)&=(-i)^mf(0,0)+\sum_{p=1}^m \Phi_p(-i)^{m-p}\\
      &=      f(0,0)(-i)^m+\sum_{p=1}^m(-i)^{m-p}\sqrt{2}CA^{-1}(\sqrt{2}A^{-1}-iI_N)^{p-1}C^*\\
      &=      f(0,0)(-i)^m+(-i)^m\sum_{p=1}^m\sqrt{2}CA^{-1}(\sqrt{2}A^{-1}-iI_N)^{-1}(\sqrt{2}A^{-1}-iI_N)^{p}(-i)^pC^*\\
      &=      f(0,0)(-i)^m+(-i)^m\sum_{p=1}^m\sqrt{2}CA^{-1}(\sqrt{2}A^{-1}-iI_N)^{-1}(i\sqrt{2}A^{-1}+I_N)^{p}C^*\\
      &=      f(0,0)(-i)^m+\\
      &\hspace{-1cm}+(-i)^m\sqrt{2}CA^{-1}(\sqrt{2}A^{-1}-iI_N)^{-1}(i\sqrt{2}A^{-1}+I_N)(I_N-(i\sqrt{2}A^{-1}+I_N)^m)(I_N-(i\sqrt{2}A^{-1}+I_N)^{-1}C^*\\
      &=      f(0,0)(-i)^m-(-i)^mC(I_N-(i\sqrt{2}A^{-1}+I_N)^m)C^*,
      \end{split}
    \]
    and so we get \eqref{newform54321} since $f(0,0)=CC^*$. To see the latter, put ${\lambda}=w=0$ in
\[
k_f({\lambda},{\nu})=\frac{\Phi({\lambda})+\Phi({\nu})^*}{1+i{\lambda}-i\overline{\nu}-{\lambda}\overline{\nu}}
\]
\end{proof}

\begin{remark}{\rm
Let $A_1=\sqrt{2}A^{-1}-iI_N$. Then $iI_N+A_1=\sqrt{2}A$ and the unique solution of the equation
\[
I_N+iA_1-iA_2-A_1A_2=0
\]
is
\[
A_2=(I_N+iA_1)(iI_N+A_1)^{-1}=(I_N+i\sqrt{2}A^{-1}+I_N)(\sqrt{2}A)^{-1}=\sqrt{2}A+iI_N=A_1^*.
\]
So we have
\begin{equation}
f(m,n)=CA_1^mA_1^{*n}C^*
  \end{equation}
and
  \begin{equation}
k_f({\lambda},{\nu})=C(I_N-{\lambda}A_1)^{-1}(I_N-\nu A_1)^{-*}C^*.
\label{bordeaux}
\end{equation}
}
\end{remark}

An immediate consequence of the previous analysis is the following:

\begin{corollary}
  Let $A=SDS^*$ where $D$ is diagonal, $D={\rm diag}\,(e^{it_1},\ldots, e^{it_N})$.
  Then,
  \begin{equation}
f(m,n)=\sum_{j=1}^Nc_j(\sqrt{2}e^{-it_j}-i)^m(\sqrt{2}e^{it_j}+i)^n
\end{equation}
for some positive numbers $c_1,\ldots, c_N$.
\end{corollary}

  \begin{example}
    \label{blabla}
Let
\[
  a=|a|e^{i\theta_a}\quad{\rm and}\quad b_a(z)=\frac{\lambda-a}{1-\lambda\overline{a}}e^{-i\theta_a}.
\]
Then,
\begin{equation}
  \label{yuiop}
  \varphi_a(z)=\frac{1+b_a(\lambda)}{1-b_a(\lambda)}=
\frac{1-|a|}{1+|a|}\frac{e^{i\theta_a}+\lambda}{e^{i\theta_a}-\lambda}.
\end{equation}

    In the case of $\varphi$ given by \eqref{yuiop} we have
    \begin{equation}
      f(m,n)=\frac{1-|a|}{1+|a|}(\sqrt{2}e^{-i\theta_a}-i)^m(\sqrt{2}e^{i\theta_a}+i)^n,\quad m,n\in \mathbb N_0,
    \end{equation}
    and in fact $m,n$ can be taken in $\mathbb Z$.
   \end{example}

   The previous arguments also hold when $H$ is not necessarily positive (but still invertible and Hermitian). The operator $A$ is now unitary with respect to the metric defined by $H$, and in
   particular may have spectrum outside the unit circle. When $-i\sqrt{2}$ is in the spectrum of $A$ the function cannot be extended to all of $\mathbb Z^2$.

\subsection{General case}

\begin{theorem}
  Assume that $\Phi_L$ is rational and analytic at infinity, with minimal realization $\Phi_L(\lambda)=D_0+C_0(\lambda I_N-A_0)^{-1}B_0$.
  Assume furthermore that $A_0$ is invertible and $\Phi_L(i)=0$. Then,
  \begin{equation}
    \label{fm0-12}
          f(m,0)=-C_0(I_N+iA_0)^{-1}A_0^{-m-1}B_0,\quad m=0,1,\ldots
          \end{equation}
        \end{theorem}

                \begin{proof}
      Assume $A_0$ invertible. We have
      \[
        \begin{split}
          \Phi_L(\lambda)&=D_0-CA_0^{-1}B_0-\lambda C_0A_0^{-1}(I_N-\lambda A_0^{-1})^{-1}A_0^{-1}B_0\\
          &=D_0-C_0A_0^{-1}B_0-\sum_{n=1}^\infty \lambda^nC_0A_0^{-n-1}B_0,
          \end{split}
        \]
where the convergence is in a neighborhood of the origin. From \eqref{q1} with $X=0$ we have
\[
  \begin{split}
    \sum_{m=0}^\infty \lambda^mf(m,0)-\frac{f(0,0)}{2}&=\frac{\Phi_L({\lambda})}{1+i\lambda}\\
    &=\left(\sum_{m=0}^\infty\lambda^m\Phi_m\right)\left(\sum_{u=0}^\infty(-i)^u\lambda^u\right).
      \end{split}
          \]
We get $\frac{f(0,0)}{2}=\Phi_L(0)=D_0-C_0A_0^{-1}B_0$ and
        \[
     f(m,0)=(D_0-C_0A_0^{-1}B_0)(-i)^m-C_0\left(\sum_{u=1}^m(-i)^{m-u}A_0^{-u-1}\right)B_0,\quad m=1,2,\ldots
          \]
          But
          \[
            \begin{split}
              \sum_{u=1}^m(-i)^{m-u}A_0^{-u-1}&=(-i)^mA_0^{-1}\sum_{u=1}^m\left(\frac{A_0^{-1}}{-i}\right)^u\\
              &=(-i)^mA_0^{-1}\frac{A_0^{-1}}{-i}\left(I_N-\left(\frac{A_0^{-1}}{-i}\right)^m\right)\left(I_N-\frac{A_0^{-1}}{-i}\right)^{-1}\\
              &=\left(-iI_N-A_0^{-1}\right)^{-1}A_0^{-2}\left((-i)^mI_N-A_0^{-m}\right)\\
               &=-A_0^{-1}(iA_0+I_N)^{-1}\left((-i)^mI_N-A_0^{-m}\right).
            \end{split}
          \]
          Hence
 \[
\begin{split}
  f(m,0)  &=\\
  &\hspace{-1cm}  =(D_0-C_0A_0^{-1}B_0)(-i)^m+C_0A_0^{-1}(I_N+iA_0)^{-1}\left((-i)^mI_N-A_0^{-m}\right)B_0\\
  &\hspace{-1cm}=(-i)^m\left(D_0-C_0A_0^{-1}B_0+C_0A_0^{-1}(I_N+iA_0)^{-1}B_0\right)-C_0(I_N+iA_0)^{-1}A_0^{-m-1}B_0\\
  &\hspace{-1cm}=(-i)^m(D_0-iC_0(I_N+iA_0)^{-1}B_0)-C_0(I_N+iA_0)^{-1}A_0^{-m-1}B_0,
          \end{split}
          \]
          since
          \[
-C_0A_0^{-1}B_0+C_0A_0^{-1}(I_N+iA_0)^{-1}B_0=-iC_0(I_N+iA_0)^{-1}B_0.
\]
The result follows since $D_0-iC_0(I_N+iA_0)^{-1}B_0)-C_0(I_N+iA_0)^{-1}A_0^{-m-1}B_0=\Phi_L(i)$.
\end{proof}

\begin{remark}{\rm
    The triple $(A_0^{-1}, A_0^{-1}B_0, C_0(I_N+iA_0)^{-1})$ is minimal  (see Lemma \ref{lemma31}), and one can recover $(A_0,B_0,C_0)$ from the values $f(m,0),$ $m=0,\ldots, N-1$.
  }
  \label{followlem}
\end{remark}

We now express $f(m,0)$ in terms of a realization of the left characteristic function $\varphi_L$ (see Definition \ref{chara123} for the latter), and begin with a preliminary computation
(see also the formulas in \cite[Theorem 1.9, p. 35]{bgk1}). We give the details for completeness. In the statement recall that $\Phi_L$ is called the left generating function (see Definition
\ref{left-right-gen}).

\begin{proposition} Let $\varphi_L(\lambda)=D+C(zI_N-A)^{-1}B$ be a minimal realization of the left characteristic function $\varphi_L$, and assume $\sqrt{2}I_N-iA$ invertible. Then,
  \begin{equation}
    \label{infty-123}
\varphi\left(\frac{\sqrt{2}{{\lambda}}}{1+i{\lambda}}\right)=D_0+C_0({\lambda}I_N-A_0)^{-1}B_0
\end{equation}
with
\begin{eqnarray}
  \label{op-1}
  A_0&=&A(\sqrt{2}I_N-iA)^{-1}\\
  B_0&=&\sqrt{2}(\sqrt{2}I_N-iA)^{-1}B\\
  C_0&=&C(\sqrt{2}I_N-iA)^{-1}\\
  D_0&=&D+iC(\sqrt{2}I_N-iA)^{-1}B,
         \label{op-4}
\end{eqnarray}
is a minimal realization of the left generating function $\Phi_L(\lambda)$.
\end{proposition}
\begin{proof}
  We can write:
\[
  \begin{split}
    \Phi({\lambda})&=\varphi\left(\frac{\sqrt{2}{{\lambda}}}{1+i{\lambda}}\right)\\
    &=D+(1+i{\lambda})C(\sqrt{2}{\lambda}I_N-(1+i{\lambda})A)^{-1}B\\
    &=D+(1+i{\lambda})C\left({\lambda}I_N-A(\sqrt{2}I_N-iA)^{-1}\right)(\sqrt{2}I_N-iA)^{-1}B\\
    &=\underbrace{D+iC(\sqrt{2}I_N-iA)^{-1}B}_{\mbox{\rm value at $\infty$}}+\\
    &\hspace{5mm}+C\left({\lambda}I_N-A(\sqrt{2}I_N-iA)^{-1}\right)(\sqrt{2}I_N-iA)^{-1}\left((1+i{\lambda})I_N-i({\lambda}I_N-A(\sqrt{2}I_N-iA)^{-1})\right)B\\
    &=D+iC(\sqrt{2}I_N-iA)^{-1}B+\\
    &\hspace{5mm}+C\left({\lambda}I_N-A(\sqrt{2}I_N-iA)^{-1}\right)(\sqrt{2}I_N-iA)^{-1}\left(I_N+iA(\sqrt{2}I_N-iA)^{-1}\right)B\\
    &=D_0+C_0({\lambda}I_N-A_0)^{-1}B_0
 \end{split}
\]
with $A_0,B_0,C_0$ and $D_0$ as in the statement of the proposition.
\end{proof}

\begin{remark}{\rm From \eqref{infty-123} we have $\Phi_L(i)=\varphi_L(\infty)=D$, which can be also checked using \eqref{op-1}-\eqref{op-4} as follows. We have
\[
  (I_N+iA_0)=I_N+iA(\sqrt{2}I_N-iA)^{-1}=\sqrt{2}(\sqrt{2}I_N-iA)^{-1}
\]
and hence,
\[
\begin{split}
  D_0-iC_0(I_N+iA_0)^{-1}B_0
  &=D+iC(\sqrt{2}I_N-iA)^{-1}B-\\
  &\hspace{5mm}-iC(\sqrt{2}I_N-iA)^{-1}\frac{1}{\sqrt{2}}(\sqrt{2}I_N-iA)\sqrt{2}(\sqrt{2}I_N-iA)^{-1}B\\
 &=D.
\end{split}
\]
}\end{remark}

      \begin{proposition} Let $\varphi_L(\lambda)=D+C(\lambda I_N-A)^{-1}B$ be a minimal realization of the left characteristic function. We have
        \begin{equation}
f(m,0)=-C(\sqrt{2}A^{-1}-iI_N)^{m}A^{-1}B,\quad m=0,1,2,\ldots
\end{equation}
\label{5-12}
        \end{proposition}

      \begin{proof}
        Using \eqref{fm0-12} and the formulas \eqref{op-1}-\eqref{op-4} for $A_0,B_0,C_0$ and $D_0$ we can write
          \[
            \begin{split}
              f(m,0)&=-C_0(I_N+iA_0)^{-1}A_0^{-m-1}B_0\\
              &=-C(\sqrt{2}I_N-iA)^{-1}(I_N+iA(\sqrt{2}I_N-iA)^{-1})^{-1}(A(\sqrt{2}I_N-iA)^{-1})^{-m-1}\times\\
              &\hspace{5mm}\times\sqrt{2}(\sqrt{2}I_N-iA)^{-1}B\\
              &=-C(\sqrt{2}I_N-iA)^{m-1}A^{-m-1}B\\
              &=-C(\sqrt{2}A^{-1}-iI_N)^{m}A^{-1}B.
            \end{split}
            \]
          \end{proof}

          To conclude we relate the moments of $\varphi_L$ and $\Phi_L$. The result is a direct consequence of the Cayley-Hamilton theorem.

          \begin{proposition}
            \label{propfin}
  In the previous notations, there exist complex numbers $a_{nk}$, $n=0,1,\ldots$ and $k=0,\ldots N-1$, such that
\[
C_0A_0^nB_0=\sum_{k=0}^{N-1}a_{nk}CA^{n+k}B,\quad n=0,1,\ldots
  \]
\end{proposition}

\begin{proof}
  We have
\[
  C_0A_0^nB_0=\sqrt{2}CA^n(\sqrt{2}I_N-iA)^{-(n+1)}B
\]
and by the Cayley-Hamilton theorem,
    \begin{equation}
     \sqrt{2}(\sqrt{2}I_N-iA)^{-(n+1)}=\sum_{k=0}^{N-1}a_{nk}A^k
      \end{equation}
    with $a_{nk}\in\mathbb C$.
      \end{proof}

      \begin{remark}
\label{remfin}
        {\rm
      One can be a bit more precise on the coefficients $a_{nk}$. Let $N_0$ be the degree of the (monic) minimal polynomial $p(\lambda)$ of $A$ and let $b_0,\ldots, b_{N_0-1}$ its coefficients:
      $p(\lambda)=\lambda^{N_0}+b_{N_0-1}\lambda^{N_0-1}+\cdots +b_0$. Write
      \[
(\sqrt{2}I_N-A)^{-1}=\sum_{k=0}^{N_0-1}c_kA^k.
        \]
        We have
        \[
          \begin{split}
            I_N&=(\sqrt{2}I_N-A)\left(\sum_{k=0}^{N_0-1}c_kA^k\right)\\
            &=\sqrt{2}\left(\sum_{k=0}^{N_0-1}c_kA^k\right)-\left(\sum_{k=1}^{N_0-1}c_{k-1}A^k\right)-c_{N_0-1}\left(-\sum_{k=0}^{N_0-1}b_kA^k\right)\\
            &=(\sqrt{2}c_0+c_{N_0-1}b_0)I_N+\sum_{k=1}^{N_0-1}(\sqrt{2}c_k-c_{k-1}+c_{N_0-1}b_k)A^k
            \end{split}
          \]
          from which we get the coefficients $c_0,\ldots, c_{N_0-1}$ via the system of equations
          \[
          \begin{split}
            \sqrt{2}c_0+c_{N_0-1}b_0&=1\\
            \sqrt{2}c_k-c_{k-1}+c_{N_0-1}b_k&=0,\quad k=1,\ldots N_0-1,
            \end{split}
          \]

          using the uniqueness of the monic minimal polynomial.
        }
      \end{remark}

      \begin{remark}{\rm In view of the previous result we reconstruct the function from the first $N$ ``moments''.}
        \end{remark}

  \section*{Acknowledgments}
  Daniel Alpay thanks the Foster G. and Mary McGaw Professorship in Mathematical Sciences, which supported this research. We thank Dr I. Paiva for the authorization to include
  Proposition \ref{ap-prop}.

\bibliographystyle{plain}
\def\cprime{$'$} \def\cprime{$'$} \def\cprime{$'$}
  \def\lfhook#1{\setbox0=\hbox{#1}{\ooalign{\hidewidth
  \lower1.5ex\hbox{'}\hidewidth\crcr\unhbox0}}} \def\cprime{$'$}
  \def\cprime{$'$} \def\cprime{$'$} \def\cprime{$'$} \def\cprime{$'$}
  \def\cprime{$'$}

\end{document}